\documentclass[11pt,a4paper]{article}

\usepackage{amssymb,amsmath}
\usepackage[top=1in,bottom=1in,left=1in,right=1in]{geometry}
\usepackage{soul}
\usepackage{graphicx,array}
\usepackage{subfigure,tabularx}
\usepackage{booktabs,caption}
\usepackage{etoolbox,verbatim,color}
\usepackage{url}
\usepackage{amsmath,amsthm,amsfonts,amssymb,amscd}
\usepackage{multirow}

\usepackage{mathtools}

\usepackage{enumitem}

\usepackage[linesnumbered,ruled]{algorithm2e}
\usepackage{appendix}
\newtheorem{theorem}{Theorem}[section]

\newtheorem{assumption}[theorem]{Assumption}
\newtheorem{corollary}[theorem]{Corollary}

\newtheorem{lemma}[theorem]{Lemma}

\theoremstyle{definition}
\newtheorem{definition}[theorem]{Definition}

\theoremstyle{definition}
\newtheorem{remark}[theorem]{Remark}

\theoremstyle{definition}
\newtheorem{example}[theorem]{Example}

\newcommand{\la}{\left\langle}
\newcommand{\ra}{\right\rangle}
\newcommand{\E}{{\mathbb{E}}}

\newcommand{\tnabla}{\widetilde{\nabla} }

\newcommand{\dom}{\operatorname{dom\,}}

\newcommand{\tnablasaga}{\widetilde{\nabla}_{\textnormal{\tiny SAGA}}}
\newcommand{\tnablasvrg}{\widetilde{\nabla}_{\textnormal{\tiny SVRG}}}
\newcommand{\tnablasarah}{\widetilde{\nabla}_{\textnormal{\tiny SARAH}}}

\newcommand{\tnablasag}{\widetilde{\nabla}_{\textnormal{\tiny SAG}}}

\newcommand{\argmin}{\ensuremath{\operatorname*{argmin}}}

\newcommand{\dist}{\ensuremath{\operatorname{dist}}}

\def\NN{\mathbb{N}}

\def\RR{\mathbb{R}}
\def\mcO{{\mathcal{O}}}
\def\mbE{{\mathbb{E}}}
\def\tx{\tilde{x}}
\newcommand{\scal}[2]{\left\langle{#1}\,,\,{#2}\right\rangle}
\newcommand{\Bscal}[2]{\Big\langle{#1},{#2} \Big\rangle}

\newcommand\numberthis{\addtocounter{equation}{1}\tag{\theequation}}

\providecommand{\keywords}[1]
{
  \small	
  \textbf{\textit{Keywords---}} #1
}

\title{Stochastic Variance-Reduced Majorization-Minimization Algorithms}
\date{May 11, 2023}
\author{
Duy-Nhat Phan \thanks{Mathematics and Statistics, UMass Lowell, MA 01854, USA. E-mail: \texttt{duynhat\_phan@uml.edu}},\quad
Sedi Bartz \thanks{Mathematics and Statistics, UMass Lowell, MA 01854, USA. E-mail: \texttt{sedi\_bartz@uml.edu}},\quad
Nilabja Guha \thanks{Mathematics and Statistics, UMass Lowell, MA 01854, USA. E-mail: \texttt{nilabja\_guha@uml.edu}} \quad
and \quad Hung M. Phan \thanks{Mathematics and Statistics, UMass Lowell, MA 01854, USA. E-mail: \texttt{hung\_phan@uml.edu}}\\[3mm]
}

\begin{document}
\maketitle
\begin{abstract}
We study a class of nonconvex nonsmooth optimization problems in which the objective is a sum of two functions: One function is the average of a large number of differentiable functions, while the other function is proper, lower semicontinuous and has a surrogate function that satisfies standard assumptions. Such problems arise in machine learning and regularized empirical risk minimization applications. However, nonconvexity and the large-sum structure are challenging for the design of new algorithms. Consequently, effective algorithms for such scenarios are scarce. We introduce and study three stochastic variance-reduced majorization-minimization (MM) algorithms, combining the general MM principle with new variance-reduced techniques. We provide almost surely subsequential convergence of the generated sequence to a stationary point. We further show that our algorithms possess the best-known complexity bounds in terms of gradient evaluations. We demonstrate the effectiveness of our algorithms on sparse binary classification problems, sparse multi-class logistic regressions, and neural networks by employing several widely-used and publicly available data sets.
\end{abstract} \hspace{10pt}


\section{Introduction} 

We focus on a class of nonsmooth and nonconvex problems of the form
\begin{equation}\label{model}
\min_{x\in\mathbb R^d}\biggl\{ F(x) := f(x) + r(x) \biggr\},
\end{equation}
where $r:\mathbb R^d \to \mathbb R\cup \{+\infty\}$ is a proper and lower semicontinuous function, and $f$ has a large-sum structure, that is,
\begin{equation}\label{eq:sum}
    f(x) = \frac{1}{n}\sum_{i=1}^nf_i(x)
\end{equation}
where $f_i$ is differentiable (possibly nonconvex). The large-sum structure captures, in particular,  \emph{regularized empirical risk}, where $f_i$ represents a loss function on a single data point and $r$ is often a nonsmooth (possibly nonconvex) function that regularizes the promotion of sparse solutions, such as $\ell_1$-norm, Geman \cite{geman1995nonlinear}, MCP \cite{zhang2010nearly}, log-sum penalty \cite{Candes_reweighting}, and exponential concave penalty \cite{Bradley1998}. Thus, problem~\eqref{model} models a broad range of optimization problems from convex (i.e., $f_i$ and $r$ are convex functions), such as logistic regression, to fully nonconvex problems (i.e., both $f_i$ and $r$ are nonconvex) such as optimizing deep neural networks. Since nonconvex optimization became indispensable in recent advances in machine learning models, we focus our attention on the fully nonconvex scenario in problem~\eqref{model}. Specifically, we are interested in the case where the number of components $n$ is extremely large since it is a key challenge in the era of big data applications.

\subsection{Motivation and Related Work}
In the convex setting, a standard method for solving the non-composite form ($r=0$) of problem \eqref{model} is the gradient descent method (\texttt{GD}). Given an initial point $x^0\in\mathbb R^d$, the iterative step of the \texttt{GD} method computes $x^{k+1}$ by
\[
x^{k+1} := x^k - \eta_k \nabla f(x^k),
\]
where $\eta_k>0$ is a stepsize. In~\eqref{eq:sum}, if the number of components $n$ is very large, each iteration of the \texttt{GD} method becomes extremely expensive since it requires the computation of the gradient for all of the components $f_i$. An effective alternative is the standard stochastic gradient method (\texttt{SGD}) \cite{robbins1951stochastic}. In this case, in each iteration, the \texttt{SGD} draws randomly $i_k$ from $[n] := \{1,2,\cdots,n\}$, and updates $x^{k+1}$ by
\[
x^{k+1} := x^k - \eta_k \nabla f_{i_k}(x^k).
\]
The advantage of the \texttt{SGD} method is that in each iteration, it only evaluates the gradient of a single component function. Consequently, the computational cost per iteration is only $1/n$ of that of the full step in the \texttt{GD} method. However, due to the {\em variance}, inadvertently generated by random sampling, the \texttt{SGD} method converges much slower than the full \texttt{GD} method. Fortunately, we can overcome this drawback by variance-reduction techniques, utilizing information regarding the gradient from previous iterations to construct a better estimation of the gradient at the current step. To date, some of the most widely applied variance reduction methods in the literature are the {\em stochastic average gradient algorithm} (\texttt{SAGA}) \cite{defazio2014saga}, the {\em stochastic variance reduced gradient} (\texttt{SVRG}) \cite{johnson2013accelerating}, and the {\em stochastic average gradient} (\texttt{SAG}) \cite{schmidt2017minimizing}. We note in passing that \texttt{SAGA} is an unbiased version of \texttt{SAG}.
Variance reduction methods inherit the advantage of low iteration cost of the \texttt{SGD} method while providing similar convergence rates of the full \texttt{GD} method in convex settings. 

Thus far, however, only several variance reduction methods have been developed in order to deal with nonconvex optimization problems possessing the large-sum structure. Furthermore, these methods mainly focus on special cases of \eqref{model}, where $r=0$, such as \cite{allen2018natasha,allen2018neon2,nguyen2019optimal}, or where $r$ is convex, such as \cite{j2016proximal,li2018simple}. For the fully nonconvex problem with an extremely large value of $n$ (such as we study here), such developments become even more challenging. Consequently, research in this direction is sparse. Several recent studies, such as \cite{Mairal2015,Lethi2017,le2020stochastic,le2022stochastic} promote stochastic methods based on difference-of-convex (DC) algorithm, developed in \cite{letthe,phacon}, or majorization minimization (\texttt{MM}), developed in \cite{lange2000optimization}. In particular, if $f_i$ is $L$-smooth and $r$ has a DC structure, that is, $r = r_1 - r_2$ with $r_1$ being proper lower semicontinuous convex and $r_2$ being convex, {\color{black} problem~\eqref{model}} can be reformulated as a DC program
\begin{equation}\label{DCprogram}
    \min_{x\in\mathbb R^d}G(x) - H(x)
\end{equation}
where $G(x) = \frac{\mu}{2}\|x\|^2 + r_1(x)$ and $H(x) = \frac{\mu}{2}\|x\|^2 - f(x) + r_2(x)$ are convex functions with $\mu\geq L$. The classic DC algorithm (\texttt{DCA}) linearizes the function $H$ iteratively and updates $x^{k+1}$ by 
\begin{equation}\label{sub-DC}
\begin{aligned}
    x^{k+1}&=\argmin_{x\in\mathbb R^d}\frac{\mu}{2}\|x\|^2 + r_1(x) - \la \mu x^k - \nabla f(x^k) + y^k, x \ra,\\
    &\color{black}=\argmin_{x\in\mathbb R^d}\frac{\mu}{2}\|x-x^k\|^2 + \la \nabla f(x^k) , x \ra
    + r_1(x) -\la y^k,x \ra\\
    &=\argmin_{x\in\mathbb R^d}\frac{\mu}{2}\|x-x^k\|^2 + \frac{1}{n}\sum_{i=1}^n
\scal{\nabla f_i(x^k)}{x}
    + r_1(x) -\scal{y^k}{x}.
\end{aligned}
\end{equation}
for some $y^k\in\partial r_2(x^k)$.

For convenience, before we continue our discussion, we set the following notations. Throughout, for $k\in\NN$, we denote by $I_k$ the index batch which is a {\em list} of (possibly repeated) indices $(i_1,i_2,\ldots,i_b)$ of fixed size $b$ where each index $i_j$ is independently and randomly chosen from $[n]$. We refer to $I_k$ as a random batch of size $b$. Let $x^0\in\mathbb{R}^d$. With a sequence $x^k\in\mathbb{R}^d$ and a sequence of batches $I_k$, we associate $x_i^{-1}=x^0$ for all $i\in[n]$, and inductively, having defined $x_i^{k-1}$, we set
\begin{equation}\label{e:xki}
x_i^{k}=\begin{cases}
x^{k},& i\in I_k,\\
x_i^{k-1},& \text{otherwise}.
\end{cases}
\end{equation}
In other words, $x_i^{k}$ is updated to $x^{k}$ if and only if $i\in I_k$.

In \cite{Lethi2017,le2020stochastic}, a stochastic version of \texttt{DCA}, named \texttt{SDCA}, was studied based
on the idea of incrementally linearizing the components $\frac{\mu}{2}\|x\|^2 - f_i(x) + r_2(x)$ of $H$. More specifically, $\mu x^k - \nabla f(x^k) + y^k$ is replaced by {\color{black} the so-called \texttt{SAG} estimator}
\begin{equation}\label{e:sag-est}
\tnablasag H(x^k) := \frac{1}{n}
    \sum_{i=1}^n\left[\mu x_i^k - \nabla f_i(x_i^k) + y_i^k\right],
    \quad\text{where~} y_i^k\in\partial r_2(x_i^k).
\end{equation}
Consequently, subproblem \eqref{sub-DC} is replaced by
\begin{equation}\label{e:sub-DC3}
x^{k+1}=\argmin_{x\in\mathbb R^d}\frac{1}{n}\sum_{i=1}^n
\Big[\frac{\mu}{2}\|x-x_i^k\|^2 + \la \nabla f_i(x_i^k) , x \ra
    + r_1(x) -\la y_i^k,x \ra\Big]
\end{equation}
where $y_i^k\in\partial r_2(x_i^k)$ for all $i$.

Recently, Le Thi et al. \cite{le2022stochastic} developed stochastic DC algorithms, named \texttt{DCA-SAGA} and \texttt{DCA-SVRG}, based on {\color{black} the so-called} \texttt{SAGA} and \texttt{SVRG} estimators for the problem in which $r$ is a DC function. {\color{black} Specifically}, \texttt{DCA-SAGA} applied to \eqref{DCprogram} successively replaces $\mu x^k - \nabla f(x^k)$ in \texttt{DCA}'s subproblem \eqref{sub-DC} by the \texttt{SAGA} stochastic gradient estimate:
\[
\begin{aligned}
\tnablasaga(\frac{\mu}{2}\|\cdot\|^2 - f)(x^k) &:= \frac{1}{b}\sum_{i\in I_k}\Big[\mu x^k - \nabla f_i(x^k) - \mu x_i^{k-1} + \nabla f_i(x_i^{k-1})\Big]\\
       &~ + \frac{1}{n}
    \sum_{i=1}^n\Big[\mu x_i^{k-1} - \nabla f_i(x_i^{k-1})\Big],
\end{aligned}
\]
which, when combined with \eqref{sub-DC}, implies that
\begin{equation}\label{e:DCA-SAGA}
\begin{aligned}
x^{k+1}=\argmin_{x\in\RR^d}~&
\frac{\mu}{2}\Big\|x-\Big(\frac{1}{b}\sum_{i\in I_k}[x^k -x_i^{k-1}]+\frac{1}{n}\sum_{i=1}^n x_i^{k-1}\Big)\Big\|^2\\
&+\frac{1}{b}\sum_{i\in I_k}\scal{\nabla f_i(x^k)-\nabla f_i(x_i^{k-1})}{x}\\
&+ \frac{1}{n}\sum_{i=1}^n\scal{\nabla f_i(x_i^{k-1})}{x} + r_1(x)-\scal{y^k}{x}.
\end{aligned}
\end{equation}
In comparison, \texttt{DCA-SVRG} replaces $\mu x^k - \nabla f(x^k)$ in \texttt{DCA}'s subproblem by the \texttt{SVRG} stochastic gradient estimate:
\begin{equation}\label{e:DCA-SVRG_est}
\begin{split}
        \tnablasvrg (\frac{\mu}{2}\|\cdot\|^2 - f)(x^k): & = \frac{1}{b}\sum_{i\in I_k}\Big[\mu x^k - \nabla f_i(x^k) - \mu\tilde x^k + \nabla f_i(\tilde x^k)\Big]
        + \mu \tilde x^k - \nabla f(\tilde x^k)\\
        & = \mu x^k - \tnablasvrg f(x^k)
    \end{split}
\end{equation}
where $\tilde x^k = x^k$ if $k\in m\NN$, and $\tilde x^{k-1}$ otherwise, where $m$ is a fixed positive integer. Consequently, the corresponding subproblem for \texttt{DCA-SVRG} is
\begin{equation}\label{e:dca-svrg-subp}
\begin{aligned}
x^{k+1}=\argmin_{x\in\RR^d}\frac{\mu}{2}\|x-x^k\|^2
&+\frac{1}{b}\sum_{i\in I_k}\scal{\nabla f_i(x^k)-\nabla f_i(\tx^k)}{x}\\
&+\scal{\nabla f(\tx^k)}{x}
+r_1(x)-\scal{y^k}{x}.
\end{aligned}
\end{equation}

In general, problem \eqref{model} can be solved by an \texttt{MM} principle such that at each iteration, a complex objective function is approximated by an upper bound which is created around the current iteration and which can be minimized effectively. This step is called the {\em majorization} step. The minimum of this upper bound (the minimization step) is then used to sequentially create another, hopefully, tighter, upper bound (another majorization step) to be minimized. More specifically, the \texttt{MM} scheme applied to \eqref{model} computes $x^{k+1}$ by
\begin{equation}\label{rule-mm}
    x^{k+1} \in\argmin_{x\in\RR^d} \frac{\mu}{2}\|x-x^k\|^2 +\la\nabla f(x^k), x-x^k\ra + u(x,x^k),
\end{equation}
where $u(x,x^k)$ is an upper bound (or surrogate, see Definition \ref{def:surrogate}) of $r$ at $x^k$.
Indeed, various deterministic approaches can  be interpreted from the \texttt{MM} point of view such as proximal or gradient-based methods \cite{combettes2011proximal,kappro,attcon,rocmon,marbre,becafa,hale2008fixed,nesterov2013gradient}, and expectation-maximization algorithms in statistics \cite{dempster1977maximum,neal1998view}. To date, many extensions of \texttt{MM} have been developed, e.g., \cite{razaviyayn2013unified,parizi2019generalized,khanh2022block,hien2023inertial,hien2022inertial,chouzenoux2022sabrina,chouzenoux2017stochastic}. However, only a few algorithms have been applied in the large-sum structure settings. In particular, Mairal~\cite{Mairal2015} introduced the {\em minimization of incremental surrogate} (\texttt{MISO}) that applies to problem \eqref{model} and successively updates $x^{k+1}$ by
\begin{equation}\label{e:miso-subp}
    x^{k+1}=\argmin \frac{1}{n}\sum_{i=1}^n\Big[\frac{\mu}{2}\|x-x_i^k\|^2 +\la\nabla f_i(x_i^k), x-x_i^k\ra + f_i(x_i^k) + u(x,x_i^k)\Big],
\end{equation}
where $u(\cdot, x_i^k)$ is a surrogate of $r$ at $x_i^k$. However, in order to study asymptotic convergence, Mairal~\cite{Mairal2015} employs a strong assumption, namely, that the approximation errors $h(x,x_i^k): = u(x,x_i^k) - r(x)$ are $L$-smooth in $x$. It is noteworthy that although \texttt{MISO} was inspired by \texttt{SAG} (see~\eqref{e:sag-est}), it does not recover \texttt{SAG}  as a special case for smooth and composite convex optimization. 

To the best of our knowledge, the incorporation of new stochastic gradient estimators \texttt{SAGA}, \texttt{SVRG}, and the {\em stochastic recursive gradient} (\texttt{SARAH}) \cite{nguyen2017sarah} into \texttt{MM} algorithms for solving the nonconvex problem \eqref{model} was not previously studied.

\subsection{Contribution and Organization}
For solving problem \eqref{model} in the case where it incorporates a large-sum structure and nonconvexity of the objective, we introduce three
{\em stochastic variance-reduced majorization-minimization} (\texttt{SVRMM}) algorithms:
\begin{itemize}[itemsep = 0pt, topsep = 0pt]
    \item \texttt{MM-SAGA} (incorporates the techniques of \texttt{SAGA});
    \item \texttt{MM-SVRG} (incorporates loop-less variants of \texttt{SVRG}); and
    \item \texttt{MM-SARAH} (incorporates loop-less variants of \texttt{SARAH}).
\end{itemize}
Unlike \texttt{MISO}, the \texttt{SVRMM} iterates on $r$ and the large-sum $f$ separately. In particular, at each iteration, \texttt{MM-SAGA}, \texttt{MM-SVRG}, and \texttt{MM-SARAH} replace the full gradient of $f$ in the deterministic \texttt{MM} by stochastic gradient estimators employing \texttt{SAGA}, loop-less \texttt{SVRG}, and loop-less \texttt{SARAH}, respectively. It is important to note that \texttt{MM-SAGA} updates the proximal term $\frac{\mu}{2}\|x-x^k\|^2$ at the current iterate $x^k$ in the same manner as in the \texttt{MM} scheme. This distinguishes \texttt{MM-SAGA} from \texttt{DCA-SAGA} when applied to the DC program \eqref{DCprogram}, where the latter's update rule consists of the proximal $\frac{\mu}{2}\|x-\bar x^k\|^2$ at $\bar x^k = \frac{1}{b}\sum_{i\in I_k}\left[x^k - x_i^{k-1}\right] + \frac{1}{n}\sum_{i=1}^nx_i^{k-1}$. In addition, we point out that \texttt{MM-SVRG} employs the loop-less \texttt{SVRG} estimator, which was shown to have superior performance \cite{kovalev2020don} when compared to the classic estimator technique \texttt{SVRG}, employed in \texttt{DCA-SVRG}. 

Under mild assumptions, we analyze the subsequential convergence for the generated sequence of the \texttt{SVRMM} algorithms. More concretely, we show that each limit point of the generated sequence is a stationary point of Problem \eqref{model}. Meanwhile, Le Thi et al. \cite{le2022stochastic} showed that each limit point $x^*$ of the generated sequence by their algorithms \texttt{DCA-SAGA} and \texttt{DCA-SVRG} is a DC critical point of $G-H$, i.e., $\partial G(x^*) \cap \partial H(x^*) \neq \emptyset$, which is weaker than the stationary point property, since $\partial F \subset \partial G - \partial H$. Furthermore, we show that our algorithms have $\mathcal O(k^{1/2})$ convergence rate with respect to the proximity to a stationary point. In order to obtain an $\epsilon$-stationary point, we show that \texttt{MM-SAGA}  and \texttt{MM-SVRG}  have complexity of $\mathcal O(n^{2/3}/\epsilon^2)$ while \texttt{MM-SARAH} has complexity of $\mathcal O(n^{1/2}/\epsilon^2)$ in terms of gradient evaluations. That is, our results are superior to those of \texttt{DCA-SAGA} and \texttt{DCA-SVRG} which have  complexity of $\mathcal O(n^{3/4}/\epsilon^2)$ and $\mathcal O(n^{2/3}/\epsilon^2)$, respectively, for finding an $\epsilon$-DC critical point. Another advantage of our methods is that we do not impose $L$-smoothness on the function $r$, it may be nonsmooth and nonconvex, but rather, in order to obtain our results, the stochastic \texttt{DCA} based algorithms require that the second DC component of $r$, namely, the component $r_2$ in the decomposition $r=r_1-r_2$, is $L$-smooth. Table \ref{complexity} contains a comparison between our new methods and \texttt{DCA-SAGA} and \texttt{DCA-SVRG} for solving nonsmooth nonconvex optimization problems, in terms of the requirement on our stepsize $1/\mu$ and our batch size $b$, in order to achieve the corresponding complexity.
{
\renewcommand{\arraystretch}{1.5}
\begin{table}[tbh!]
\centering
\caption{\texttt{SVRMM} vs stochastic-based \texttt{DCA} }\label{tab:sum1}
\begin{tabular}{@{}|l|llll|@{}} 
\hline
 Method &  Requirement  & Stepsize $1/\mu$ & Batch size & Complexity\\  \hline
\texttt{DCA-SAGA} \cite{le2022stochastic} & $\frac{n\sqrt{n+b}}{b^2}\leq \frac{1}{4}\frac{2\mu-L}{\mu+L}$ & $\frac{1}{2L}$ & $b = \lfloor 2^{1/4}2n^{3/4}\rfloor$ & $\mcO (n^{3/4}/\epsilon^2)$\\
\texttt{DCA-SVRG} \cite{le2022stochastic} & $\frac{m}{\sqrt{b}}\leq \frac{1}{4\sqrt{e-1}}\frac{2\mu-L}{\mu+L}$ & $\frac{1}{2L}$ & $b = \lfloor n^{2/3}\rfloor$, $m=\lfloor \frac{\sqrt{b}}{4\sqrt{e-1}}\rfloor$ & $\mcO(n^{2/3}/\epsilon^2)$ \\ 
\texttt{MM-SAGA} (new) & $\frac{n}{b^{3/2}}\leq \frac{1}{4}\frac{2\mu-L}{L}$ & $\frac{1}{L}$ & $b = \lfloor 4^{2/3}n^{2/3}\rfloor$ & $\mathcal O(n^{2/3}/\epsilon^2)$\\
\texttt{MM-SVRG} (new) & $\frac{m}{\sqrt{b}}\leq \frac{1}{4}\frac{2\mu-L}{L}$ & $\frac{1}{L}$& $b = \lfloor n^{2/3}\rfloor$, $m=\lfloor \frac{\sqrt{b}}{4}\rfloor$ & $\mcO(n^{2/3}/\epsilon^2)$\\
\texttt{MM-SARAH} (new) & $\frac{m}{b}\leq \frac{1}{4}\frac{(2\mu-L)^2}{L^2}$ & $\frac{1}{L}$ & $b = \lfloor n^{1/2} \rfloor$, $m=\lfloor \frac{b}{4}\rfloor$ & $\mcO(n^{1/2}/\epsilon^2)$\\
\hline
\end{tabular}\label{complexity}
\vspace{-0.1in}

\end{table}
}

In Table~\ref{tab:sum1}, the stepsizes and batchsizes of \texttt{DCA-SAGA} and \texttt{DCA-SVRG} are taken from \cite{le2022stochastic}, and the stepsizes and batchsizes for the \texttt{MM-SAGA}, \texttt{MM-SVRG}, and \texttt{MM-SARAH} are chosen in order to achieve the optimal order of complexity. We also note by passing that by setting the stepsize $\frac{1}{\mu}=\frac{1}{2L}$, we obtain the same complexity in all three \texttt{SVRMM} algorithms.

Finally, we apply our algorithms to solve three problems: sparse binary classification with nonconvex loss and regularizer, sparse multi-class logistic regression with nonconvex regularizer, and feedforward neural network training, in order to illustrate their applicability and efficiency.

The paper is organized as follows. In Section~\ref{sec:prelnnopt}, we present basic concepts and properties
in nonconvex optimization. In Section~\ref{SVRMM}, we present our algorithms. We analyze the convergence properties of our methods in Section~\ref{convergence}. In Section \ref{numerical_exp}, we provide a demonstration by numerical experiments, followed by conclusions in Section~\ref{conclusion}.

\section{Preliminaries}
 \label{sec:prelnnopt}

We follow standard notations such as in~\cite{bauschke_combettes_2017,VariationalAnalysis}. Throughout, $\scal{\cdot}{\cdot}$ denotes the inner product in $\mathbb R^d$ with induced norm $\|\cdot\|$ defined by $\|x\| = \sqrt{\la x,x\ra}, x\in\mathbb R^d$. We set $\mathbb R_+ = \{r \in\mathbb R:r\geq 0\}$. For a nonempty closed set $\mathcal C\subset\mathbb R^d$, the distance of $x$ from $\mathcal C$ is defined by $\dist(x,\mathcal C) = \inf_{y\in\mathcal C}\|x-y\|$. An extended-real-valued function $g:\mathbb R^d\to \mathbb R\cup \{+\infty\}$ is said to be proper if its domain, the set ${\rm dom}\,g = \{x\in\mathbb R^d:g(x)<+\infty\}$, is nonempty. We say that $g$ is lower semicontinuous if, at each $x^*\in\mathbb R^d$,
 \begin{equation*}
     g(x^*)\leq \liminf_{x\to x^*}g(x).
 \end{equation*}
 Let $\alpha\in\mathbb R$. We say that $g$ is $\alpha$-convex if $g-\frac{\alpha}{2}\|\cdot\|^2$ is convex, equivalently, if 
 \begin{equation*}
     g((1-\lambda)x+\lambda y) \leq (1-\lambda) g(x) + \lambda g(y) - \frac{\alpha}{2}\lambda(1-\lambda)\|x-y\|^2\quad\forall x,y\in\mathbb R^d, \lambda\in[0,1].
 \end{equation*}
 In particular, $g$ is convex if and only if $g$ is $0$-convex. In the case where the function $g$ is $\alpha$-convex, we say that $g$ is $\alpha$-strongly convex if $\alpha>0$ and we say that $g$ is $\alpha$-weakly convex if $\alpha<0$. 

\begin{definition}[Fr\'echet and limiting subdifferential]\cite[Definition 8.3]{VariationalAnalysis}
\label{def:dd}
Let $g: \mathbb R^d\to \mathbb R\cup \{+\infty\} $ be a proper lower semicontinuous function. 
\begin{enumerate}[label={\rm(\alph*)}]
\item For each $x\in{\rm dom}\,g,$ we denote by $\hat{\partial}g(x)$ 
the Fr\'echet subdifferential of $g$ at $x$. It contains all of the vectors
$v\in\mathbb R^d$ which satisfy
\[
\liminf_{y\ne x,y\to x}\frac{1}{\left\Vert y-x\right\Vert }\left(g(y)-g(x)-\left\langle v,y-x\right\rangle \right)\geq 0.
\]
If $x\not\in{\rm dom}\:g,$ we set $\hat{\partial}g(x)=\emptyset.$  
\item The limiting-subdifferential $\partial g(x)$ of $g$ at $x\in{\rm dom}\:g$
is defined by
\[
\partial g(x) := \left\{ v\in\mathbb R^d:\exists x^{k}\to x,\,g\big(x^{k}\big)\to g(x),\,v^{k}\in\hat{\partial}g\big(x^{k}\big),\,v^{k}\to v\right\} .
\]
\end{enumerate}
\end{definition}
 If $g$ is convex, then the Fr\'echet and the limiting subdifferential coincide with the convex subdifferential:
\[
\partial g (x) = \left\{v : g(y) \ge g(x) + \la y-x,v\ra, \forall y \in \mathbb R^d\right\}.
\]

\begin{definition}[$L$-Lipschitz] A mapping $T:D\subset \mathbb{R}^d\to\mathbb{R}^k$ is said to be $L$-Lipschitz, $L\geq 0$, if
\[
\forall x,y\in D,\quad
\|Tx-Ty\|\leq L\|x-y\|. 
\]
\end{definition}

\begin{definition}[$L$-smooth]
Let $g: \mathbb R^d\to \mathbb R$. We say that $g$ is $L$-smooth  if it is differentiable and its gradient, $\nabla g$, is $L$-Lipschitz.
\end{definition}

We now recall several useful basic facts.
\begin{lemma}\label{convexity-smoothness}
Let $g$ and $h$ be proper and lower semicontinuous. Then,
    \begin{enumerate}[label={\rm(\alph*)}]
    \item $0\in\partial g(\bar x)$ if $g$ attains a local minimum at $\bar x\in{\rm dom}\:g$.
    \item $\partial f(\bar x) = \partial g(\bar x) + \nabla h(\bar x)$ if $f = g+h$ and $h$ is continuously differentiable in a neighborhood of $\bar x$.
    
    \item $g(x)\geq g(y) + \rho\|x-y\|^2$ if $g$ is convex and $y$ is defined by
    \[
    y = \argmin_z\left\{g(z) + \frac{\rho}{2}\|z-x\|^2\right\}.
    \]
    \item $|g(x) - g(y) - \la \nabla g(y),x -y\ra| \leq \frac{L}{2}\|x-y\|^2$ $\forall x,y \in\mathbb R^d$ if $g$ is $L$-smooth.
    \end{enumerate}
\end{lemma}
\begin{proof}
    a) see, e.g., \cite[Theorem 8.15]{VariationalAnalysis}. b) See, e.g., \cite[Exercise 8.8]{VariationalAnalysis}. c) See, e,g., \cite[Theorem 6.39]{beck2017first}. d) See, e.g., \cite[Lemma 1.2.3]{nesterov2018lectures}.
\end{proof}
\begin{definition}[Stationary point]
\label{def:type2}
We say that $x^{*}\in \dom g$ a stationary point of $g$ if $0\in\partial g\left(x^{*}\right).$ 
\end{definition}

\begin{definition}[$\epsilon$-stationary point]
A point $x^*$ is said to be an $\epsilon$-stationary point of $g$ if 
\[
\dist\big(0,\partial g(x^*)\big)\leq \epsilon.
\]
\end{definition}
The following lemma is a fundamental tool in our convergence analysis.

\begin{lemma}[Supermartingale convergence
{\rm\cite[Theorem 1]{robbins1971convergence}}]
\label{supermartingale}
Let $\{Y_k\}, \{Z_k\}$, and $\{W_k\}$ be three sequences of random variables and let $\mathcal F_k$ be sets of random variables such that $\mathcal F_k\subset \mathcal F_{k+1}$ for all $k$. Assume that
\begin{enumerate}[label={\rm(\alph*)}]
    \item The random variables $\{Y_k\}, \{Z_k\}$, and $\{W_k\}$ are nonnegative and are functions of random variables in $\mathcal F_k$;
    \item $\E\big[Y_{k+1}|\mathcal F_k\big]\leq  Y_k - Z_k + W_k$ for each $k$;
    \item $\sum_{k=0}^{+\infty}W_k < +\infty$ with probability $1$.
\end{enumerate}
Then, $\sum_{k=0}^{+\infty}Z_k < +\infty$, and $\{Y_k\}$ converges to a nonnegative random variable, almost surely.
\end{lemma}

\section{Stochastic Variance-Reduced Majorization-Minimization}\label{SVRMM}
In this section, we introduce three {\em stochastic variance-reduced majorization-minimization} (\texttt{SVRMM}) algorithms for solving problem \eqref{model}. To this end, we define surrogate functions as follows. 
\begin{definition}[Surrogate function]\label{d:surrogate}
\label{def:surrogate}
A function $u:\mathbb R^d\times\mathbb{R}^d \to \mathbb R\cup \{+\infty\}  $ is said to be a surrogate function of $r:\mathbb R^d\to \mathbb R\cup \{+\infty\}$
if 
    \begin{enumerate}[label={\rm(\alph*)}]
    \item $u(y,y) = r(y)$ for all $y\in\mathbb R^d$;
    \item $u(x,y) \geq r(x)$ for all $x,y\in\mathbb R^d$.
    \end{enumerate}
\end{definition}
We introduce our first \texttt{SVRMM} algorithm, which we call \texttt{MM-SAGA}. It combines the deterministic \texttt{MM} and the \texttt{SAGA}-style of stochastic gradient update. In particular, we replace the full gradient $\nabla f(x^k)$ in the deterministic \texttt{MM} \eqref{rule-mm} with the stochastic gradient estimate $\tnablasaga f(x^k)$ \eqref{e:tnablasaga}.
\texttt{MM-SAGA} is formally detailed in Algorithm~\ref{mm-saga}.

\begin{algorithm}[H] 
\caption{\texttt{MM-SAGA}}
\label{mm-saga}
\SetKwInput{KwInput}{Input}
\DontPrintSemicolon

\KwInput{$x^0\in\mathbb R^d, v^0 = \nabla f(x^0)$, $\nabla f_i(x^{-1}_i) = \nabla f(x^0)$ for $i=1,\ldots,n$, $k=0$, $\mu>\frac{L}{2}$, a batch size $b$, and a surrogate function $u$ of $r$.}
\Repeat{Stopping criterion.}{
Choose a random batch $I_k$ of size $b$.

Set $x^k_i$ by \eqref{e:xki} and define the \texttt{SAGA} gradient estimate
\begin{equation}\label{e:tnablasaga}
    \tnablasaga f(x^k) = \frac{1}{b}\sum_{i\in I_k}\left(\nabla f_i(x^k) - \nabla f_i(x^{k-1}_i)\right) +  v^k.
\end{equation}

Compute
\begin{equation}\label{rule-mm-saga}
 x^{k+1} = \argmin_{x\in\mathbb R^d} \frac{\mu}{2}\|x-x^k\|^2 +\la \tnablasaga f(x^k), x\ra  + u(x,x^k).
\end{equation}

Update
\[
v^{k+1} := \frac{1}{n}\sum_{i\in I_k}\nabla f_i(x^k_i) = \frac{1}{n}\sum_{i\in I_k}\left(\nabla f_i(x^k) - \nabla f_i(x^{k-1}_i)\right) + v^k.
\]
Set $k\leftarrow k+ 1$.
}
\end{algorithm}

\begin{remark}[\texttt{MM-SAGA} vs. \texttt{DCA-SAGA}]
The update rule \eqref{e:DCA-SAGA} of \texttt{DCA-SAGA} \cite{le2022stochastic}
can be expressed as
\[
 \min_{x\in\mathbb R^d}\frac{\mu}{2}\|x - \bar x^k\|^2 +  \la\tnablasaga f(x^k),x\ra + r_1(x) - \la y^k,x\ra,
\]
where $\bar x^k = \frac{1}{b}\sum_{i\in I_k}\left(x^k - x_i^{k-1}\right) + \frac{1}{n}\sum_{i=1}^nx_i^{k-1}$ and $y^k\in\partial r_2(x^k)$. This update is different from our update rule \eqref{rule-mm-saga} in \texttt{MM-SAGA}  which employs the proximal term $\frac{\mu}{2}\|x-x^k\|^2$, in the same manner as in the deterministic \texttt{MM}.

It is worth mentioning an advantage of \texttt{MM-SAGA} over \texttt{DCA-SAGA} in terms of memory storage, which can be described as follows. For \texttt{MM-SAGA}, one employs $x_i^k$ solely for computing $\nabla f_i(x_i^k)$ only. In other words, \eqref{e:tnablasaga} and \eqref{rule-mm-saga} are completely defined and updated using only the gradient values. Thus, it suffices to keep only the value $\nabla f_i(x_i^k)$, and then discard $x_i^k$. In many problems, such as binary classifications and multi-class logistic regressions, each gradient $\nabla f_i$ is a scalar multiple of the data point $i$, where the scalar can be updated in each iteration. In such cases, one may only store the weights, instead of the full gradient $\nabla f_i$. In contrast, for \texttt{DCA-SAGA}, in addition to the values of gradients $\nabla f_i(x_i^k)$, one must store all vectors $x_i^k$ in order to compute $\bar{x}^k$ in each iteration.
\end{remark}

Our second \texttt{SVRMM} algorithm, named \texttt{MM-SVRG}, is inspired by a loop-less \texttt{SVRG} estimator. Specifically, we replace the full gradient $\nabla f(x^k)$ in the deterministic \texttt{MM} \eqref{rule-mm} by the loop-less \texttt{SVRG} stochastic gradient estimate $\tnablasvrg f(x^k)$~\eqref{update-svrg} (see \cite{kovalev2020don}). 
\texttt{MM-SVRG} is formally detailed in Algorithm \ref{mm-svrg}.

\begin{algorithm}[H] 
\caption{\texttt{MM-SVRG}}
\label{mm-svrg}
\SetKwInput{KwInput}{Input}
\DontPrintSemicolon

\KwInput{$x^0\in\mathbb R^d, \tilde x^{-1} = x^0$, positive integer $m$, a batch size $b$, $k=0$, $\mu>\frac{L}{2}$, and a surrogate function $u$ of $r$.}
\Repeat{Stopping criterion.}{
Choose a random batch $I_k$ of size $b$.

Define the (loop-less) \texttt{SVRG} gradient estimate
\begin{equation}\label{update-svrg}
    \tnablasvrg f(x^k) = \frac{1}{b}\sum_{i\in I_k}\left(\nabla f_i(x^k) - \nabla f_i(\tilde x^k)\right) +  \nabla f(\tilde x^k)
\end{equation}
where $\tilde x^k = x^k$ w.p. $1/m$ and $\tilde x^{k-1}$ otherwise.

Compute
\begin{equation}\label{rule-mm-svrg}
x^{k+1}= \argmin_{x\in\mathbb R^d} \frac{\mu}{2}\|x-x^k\|^2 +\la \tnablasvrg f(x^k), x\ra  + u(x,x^k).
\end{equation}

Set $k\leftarrow k+1$ 
}
\end{algorithm}

\begin{remark}
It is worth noting that \texttt{MM-SVRG} employs the loop-less \texttt{SVRG} estimator which was demonstrated to have practical advantages in \cite{kovalev2020don} when compared to the estimator \eqref{e:DCA-SVRG_est} used in \texttt{DCA-SVRG} \cite{le2022stochastic}. 
\end{remark}

Our third \texttt{SVRMM} algorithm is named \texttt{MM-SARAH}, in which we replace the gradient $\nabla f(x^k)$ in the deterministic \texttt{MM} \eqref{rule-mm} with a loop-less variant of \texttt{SARAH} \eqref{update-sarah}.
\texttt{MM-SARAH} is formally detailed in Algorithm~\ref{mm-sarah}.

\begin{algorithm}[ht!] 
\caption{\texttt{MM-SARAH}}
\label{mm-sarah}
\SetKwInput{KwInput}{Input}
\DontPrintSemicolon

\KwInput{$x^{-1}=x^0\in\mathbb R^d, \tnablasarah f(x^{-1}) = \nabla f(x^0)$, positive integer $m$, a batch size $b$, $k=0$, $\mu>\frac{L}{2}$, and a surrogate function $u$ of $r$.}
\Repeat{Stopping criterion.}{
Choose a random batch $I_k$ of size $b$.

Define the (loop-less) \texttt{SARAH} gradient estimate
\begin{equation}\label{update-sarah}
    \tnablasarah f(x^k) = \begin{cases}
	   \nabla f(x^k)  & \textrm{w.p. } 1/m, \\
	    \frac{1}{b} \sum\limits_{i \in I_k}\left(\nabla f_i(x^k) - \nabla f_i(x^{k-1})\right) + \tnablasarah f(x^{k-1}) & \textrm{otherwise}
    \end{cases}
\end{equation}

Compute
\begin{equation}\label{rule-mm-sarah}
x^{k+1}=\argmin_{x\in\mathbb R^d} \frac{\mu}{2}\|x-x^k\|^2 +\la \tnablasarah f(x^k), x\ra  + u(x,x^k).
\end{equation}

Set $k\leftarrow k+ 1$
}
\end{algorithm}

\bgroup
Finally, for convenience, we summarize all of the algorithms under discussion in the following three tables. First, the gradient estimators are detailed as follows.
\begin{equation*}
\arraycolsep 10pt
\def\arraystretch{2}
\begin{array}{ |l|c|}
\hline
\displaystyle
\tnablasaga f(x^k)=\frac{1}{b}{\sum_{i\in I_k}}\left(\nabla f_i(x^k)-\nabla f_i(x_i^{k-1}\right) 
+ \frac{1}{n}\sum_{i=1}^n\nabla f_i(x_i^{k-1}) & \text{see~}\eqref{e:tnablasaga}\\
\hline
~&~\\[-5ex]
\displaystyle
\begin{aligned}
\tnablasvrg f(x^k)=\frac{1}{b}\sum_{i\in I_k}\left( \nabla f_i(x^k)-\nabla f_i(\tx^k)\right)
+ \nabla f(\tx^k)\\ 
    \text{where~~}\tx^k = x^k\text{~w.p~} 1/m,\ \tx^{k-1}\text{~w.p.~} 1-1/m
\end{aligned}
& \text{see~}{\eqref{e:DCA-SVRG_est},\eqref{update-svrg}} 
\\[+1ex]
\hline
~&~\\[-4ex]
\tnablasarah f(x^k) =
\begin{cases}
\nabla f(x^k),&\text{~w.p.~} 1/m, \\
\frac{1}{b} \sum\limits_{i \in I_k}\left(\nabla f_i(x^k) - \nabla f_i(x^{k-1})\right) + \tnablasarah f(x^{k-1}),&\text{~w.p.~} 1-1/m
\end{cases}
&
\text{see~}\eqref{update-sarah} 
\\[+3ex]
\hline
\end{array}
\end{equation*}

The next table details iteration updates for algorithms from the literature.
\[
\arraycolsep 10pt
\def\arraystretch{2}
\begin{array}{ |c|c|c|}
\hline
  & x^{k+1} \text{~ is a minimizer of} & \text{see}\\
\hline
\texttt{DCA}     & 
\frac{\mu}{2}\|x-x^k\|^2 + 
\frac{1}{n}\sum_{i=1}^n
\la \nabla f_i(x^k) , x \ra
    + r_1(x) -\la y^k,x \ra,
    \quad y^k\in\partial r_2(x^k)
& \eqref{sub-DC}
    \\[+1ex]
\hline
\texttt{SDCA}     &\frac{1}{n}\sum_{i=1}^n
\Big[\frac{\mu}{2}\|x-x_i^k\|^2 + \la \nabla f_i(x_i^k) , x \ra
    + r_1(x) -\la y_i^k,x \ra\Big],
    \quad y^k_i\in\partial r_2(x^k_i) & 
    \eqref{e:sub-DC3}
    \\[+1ex]
\hline
\texttt{DCA-SAGA} & 
\begin{aligned}
\frac{\mu}{2}\Big\|x-\Big(\frac{1}{b}\sum_{i\in I_k}[x^k &-x_i^{k-1}] + \frac{1}{n}\sum_{i=1}^n x_i^{k-1}\Big)\Big\|^2 \\
&+\scal{\tnablasaga f(x^k)}{x}
+ r_1(x)-\scal{y^k}{x}
\end{aligned}
& \eqref{e:DCA-SAGA}
\\ 
~&~&\\[-4.5ex]
\hline
~&&\\[-4ex]
\texttt{DCA-SVRG} &
\frac{\mu}{2}\|x-x^k\|^2
+\scal{\tnablasvrg f(x^k)}{x}
+r_1(x)-\scal{y^k}{x}
& \eqref{e:dca-svrg-subp} \\[+3ex]
\hline
\texttt{MISO} &
\frac{1}{n}\sum_{i=1}^n
\Big[\frac{\mu}{2}\|x-x_i^k\|^2 + \scal{\nabla f_i(x_i^k)}{x}
+ u(x,x^k_i)\Big] &
\\
& \text{$u(\cdot,\cdot)$ is a surrogate of $r(\cdot)$, (see~Definition~\ref{d:surrogate}).}
& \eqref{e:miso-subp}\\
\hline
\end{array}
\]
\egroup

\bgroup
Finally, iteration updates for our three \texttt{SVRMM} algorithms are summarized in the following table, where $u(\cdot,\cdot)$ is a surrogate function for $r(\cdot)$, see~Definition~\ref{d:surrogate}.
\[
\arraycolsep 10pt
\def\arraystretch{2}
\begin{array}{ |c|c|c|}
\hline
  & x^{k+1} \text{~ is a minimizer of} & \text{see}\\
\hline
\texttt{MM-SAGA} & 
\displaystyle\frac{\mu}{2}\|x-x^k\|^2 +\la \tnablasaga f(x^k), x\ra  + u(x,x^k)
& \text{Alg.~\ref{mm-saga}}\\[+1ex]
\hline
~&&\\[-4ex]
\texttt{MM-SVRG} &
\frac{\mu}{2}\|x-x^k\|^2
+\scal{\tnablasvrg f(x^k)}{x}
+u(x,x^k)
& \text{Alg.~\ref{mm-svrg}}\\
~&&\\[-4ex]
\hline
\texttt{MM-SARAH} &
\displaystyle\frac{\mu}{2}\|x-x^k\|^2 +\la \tnablasarah f(x^k), x\ra  + u(x,x^k)
& \text{Alg.~\ref{mm-sarah}}
\\[+1ex]
\hline
\end{array}
\]
\egroup

\section{Convergence analysis}\label{convergence}
We now focus on convergence analysis of our \texttt{SVRMM} algorithms. To this end, throughout this paper, we assume the following basic assumptions regarding problem \eqref{model}.
Such assumptions are standard in optimization literature, see, e.g., \cite{nesterov2003introductory,pham2020proxsarah}. 
\begin{assumption}\label{basic:assumption}
\begin{enumerate}[label={\rm(\alph*)}]
    \item $F$ is bounded from below, that is, $F^* = \inf_{x\in\mathbb R^d}F(x) >-\infty$, and 
    \[
    {\rm dom}\:F = {\rm dom}\:f\cap {\rm dom}\:r\neq\emptyset.  
    \]
    \item $f_i$ is continuously differentiable and there exists a positive constant $L$ such that
    \begin{equation}\label{average-smoothness}
        \frac{1}{n}\sum_{i=1}^n\|\nabla f_i(x) - \nabla f_i(y)\|^2 \leq L^2\|x-y\|^2,\quad\forall x,y\in {\rm dom}\: f.
    \end{equation}
\end{enumerate}
\end{assumption}
We note that \eqref{average-smoothness} is satisfied if each component function $f_i$ is $L_i$-smooth. Indeed, in this case,
\[
\frac{1}{n}\sum_{i=1}^n\|\nabla f_i(x) - \nabla f_i(y)\|^2 \leq \frac{1}{n}\sum_{i=1}^nL_i^2\|x-y\|^2,
\]
which implies that \eqref{average-smoothness} holds with $L^2 = \frac{1}{n}\sum_{i=1}^nL_i^2$. It is also worth noting that \eqref{average-smoothness} implies $L$-smoothness of $f$. Indeed, 
\[
\|\nabla f(x) - \nabla f(y)\|^2 = \left\|\frac{1}{n}\sum_{i=1}^n[\nabla f_i(x) - \nabla f_i(y)]\right\|^2\leq \frac{1}{n}\sum_{i=1}^n\left\|\nabla f_i(x) - \nabla f_i(y)\right\|^2\leq L^2\|x-y\|^2,
\]
which implies that $\|\nabla f(x) - \nabla f(y)\|\leq L\|x-y\|$, i.e., $f$ is $L$-smooth.

Our following assumption is regarding the surrogate function $u$ of $r$, see Definition~\ref{d:surrogate}.
\begin{assumption}\label{assump:u}
\begin{enumerate}[label={\rm(\alph*)}]
\item For every $x$, $u(x,\cdot)$ is continuous in $y$.
\item For every $y$, $u(\cdot,y)$ is {\color{black} lower semicontinuous} and convex.
 \item There exists a function $\bar h:\mathbb{R}^d\times\mathbb{R}^d\to\mathbb{R}\cup 
 \{+\infty\}$ such that for every $y\in\mathbb{R}^d$, $\bar h(\cdot,y)$ is continuously differentiable at $y$ with $\nabla \bar h(\cdot,y)(y)=0$, and  
 the approximation error satisfies
\begin{equation}
\label{lemma:h_property} 
u(\cdot,y)-r(\cdot)\leq \bar h(\cdot,y).
\end{equation}
\end{enumerate}
\end{assumption}

\begin{remark}
    We note that Assumption \ref{assump:u}(c) encapsulates surrogate functions previously studied by Mairal \cite[Definition 2.2]{Mairal2015} where it is assumed that the approximation error $h(\cdot,y): = u(\cdot,y) - r(\cdot)$ is $L$-smooth and $\nabla h(\cdot,y)(y) = 0$. Indeed, in this case, we simply set $\bar{h}(\cdot,y)=u(\cdot,y)-r(\cdot)$.
\end{remark}

We recall several classes of functions $r$ which satisfy Assumption \ref{assump:u}. In particular, we recall functions with $L$-smooth approximation error $h(\cdot,y)$ such that $\nabla h(\cdot,y)(y) = 0$. Additional examples can be found in \cite{Mairal2015,khanh2022block}.
\begin{example}[Proximal surrogates]\label{ex:prox-surr} If $r$ is $\alpha$-weakly convex, it is natural to consider the surrogate
\[
u(x,y) = r(x) + \frac{-\alpha}{2}\|x - y\|^2.
\]
\end{example}

\begin{example}[Lipschitz gradient surrogates]
\label{ex:lip-grad-surr}
If $r$ is $L$-smooth, a natural surrogate is
\[
u(x,y) = r(y) + \la \nabla r(y), x - y\ra + \frac{\rho}{2}\|x-y\|^2,
\]
for some $\rho\geq L$.
\end{example}

\begin{example}[DC surrogates]\label{ex:dc-surr}
If $r$ is a DC function, that is, $r = r_1 - r_2$ where $r_2$ is $L$-smooth, we consider the surrogate function
\[
u(x,y) = r_1(x) - \Big[\la \nabla r_2(y), x - y\ra + r_2(y)\Big].
\]
We now provide an example in which the approximation error function is nonsmooth, however, Assumption \ref{assump:u} is satisfied. 
\end{example}

\begin{example}[Composite surrogates]\label{ex:comp_surr}
Consider a class of functions of the form
\[
r(x) = \sum_{i=1}^m\eta_i\big(g_i(x_i)\big),
\]
where $x$ is decomposed into $m$ blocks $x=(x_1,...,x_m)$ with $x_i\in\mathbb R^{d_i}$, $\sum_{i=1}^md_i=d$, and where $g_i:\mathbb R^{d_i}\to \mathbb R$ are convex and Lipschitz continuous with a common Lipschitz constant $L_g$, and $\eta_i:\mathbb R\to \mathbb R$ are concave and smooth with a common smoothness constant $L_\eta$ on the image $g_i(\mathbb R)$. This class includes composite functions, in particular, several existing sparsity-reduced regularizers, which are nonconvex and nonsmooth approximations of the $\ell_0$-norm or $\ell_{q,0}$-norm, see, e.g., \cite{Bradley1998,Candes_reweighting}. Since $\eta_i$ is concave, we can setup a surrogate function $u$ for $r$ as follows
\begin{equation*}
    u(x,y) = r(y) + \sum_{i=1}^m \eta'_i\big(g_i(y_i)\big)\big(g_i(x_i) - g_i(y_i)\big).
\end{equation*}
Since $\eta_i$ is $L_\eta$-smooth on the image $g_i(\RR)$, it follows from Lemma \ref{convexity-smoothness} (d) that
\begin{equation*}
    r(x) \geq r(y) + \sum_{i=1}^m\eta'_i\big(g_i(y_i)\big)\big(g_i(x_i) - g_i(y_i)\big) - \sum_{i=1}^m\frac{L_\eta}{2}\|g_i(x_i) - g_i(y_i)\|^2,
\end{equation*}
which, when combined with the $L_g$-Lipschitz continuity of $g_i$, implies that
\[
    u(x,y) - r(x) \leq \sum_{i=1}^m\frac{L_\eta}{2}\|g_i(x_i) - g_i(y_i)\|^2\leq \frac{L_\eta L_g^2}{2}\|x-y\|^2.
\]
Therefore, Assumption \ref{assump:u} (c) is satisfied when we set $\bar h(x,y) = \frac{L_\eta L_g^2}{2}\|x-y\|^2$.
\end{example}

The following lemma will be useful in our convergence analysis.
\begin{lemma}
\label{l:2023feb01}
Let $\xi_i\in\mathbb R^n$ and set $\bar{\xi}=\frac{1}{n}\sum_{i=1}^{n}\xi_i$. Suppose that $I$ is a list of length $b$ of (possibly repeated) indices $(i_1,\ldots,i_b)$ where each $i_k$ is chosen independently and randomly from $[n]$. Then
\begin{equation}
\mbE\left\|\frac{1}{b}\sum_{i\in I}\xi_i - \bar{\xi}
\right\|^2=\frac{1}{bn}\sum_{i=1}^{n}\|\xi_i\|^2-
\frac{1}{b}\|\bar{\xi}\|^2,
\end{equation}
where the expectation is taken over all possible choices of $I$.
\end{lemma}
\proof We begin with the observation that
\begin{align*}
\mbE\left\|\frac{1}{b}\sum_{i\in I}\xi_i - \bar{\xi}
\right\|^2
&=\frac{1}{b^2}\mbE\left\|\sum_{i\in I}(\xi_i - \bar{\xi})
\right\|^2\\
&=\frac{1}{b^2}\mbE\left(
\sum_{i\in I}\|\xi_i-\bar{\xi}\|^2
+
\sum_{i,j\in I,\atop i\neq j}\la\xi_i - \bar{\xi},\xi_j-\bar{\xi}\ra
\right).
\end{align*}
Since $i,j$ are chosen independently and randomly from $[n]$, then for each pair $(i,j)$ such that $i\neq j$,
\begin{equation*}
\mbE\scal{\xi_i-\bar{\xi}}{\xi_j-\bar{\xi}}
=\Bscal{\mbE(\xi_i-\bar{\xi})}{\mbE(\xi_j-\bar{\xi})}
=0
\end{equation*}
and
\begin{align*}
\mbE\|\xi_i-\bar{\xi}\|^2
&=\mbE(\|\xi_i\|^2 - 2\scal{\xi_i}{\bar{\xi}} + \|\bar{\xi}\|^2)\\
&=\mbE\|\xi_i\|^2 - 2\scal{\mbE\xi_i}{\bar{\xi}} + \|\bar{\xi}\|^2
=\mbE\|\xi_i\|^2 - \|\bar{\xi}\|^2\\
&=\frac{1}{n}\sum_{i=1}^{n}\|\xi_i\|^2-\|\bar{\xi}\|^2.
\end{align*}
By recalling that $I$ contains $b$ indices, it follows that
\begin{equation*}
\mbE\left\|\frac{1}{b}\sum_{i\in I}\xi_i - \bar{\xi}
\right\|^2=
\frac{1}{bn}\sum_{i=1}^{n}\|\xi_i\|^2-
\frac{1}{b}\|\bar{\xi}\|^2
\end{equation*}
which completes the proof.
\endproof

Before we proceed to our convergence analysis, for  simplicity, we associate with each of our algorithms a sequence of random variables $\{\Upsilon^k\}$, which is defined by
\begin{equation}\label{e:Upsilon}
\arraycolsep 10pt
\def\arraystretch{1.5}
\begin{array}{l|l}
\hline
\text{Method} & \Upsilon^k\\
\hline
\texttt{MM-SAGA}& \frac{1}{bn}\sum_{i=1}^n\left\|\nabla f_i(x^k) - \nabla f_i(x_i^{k-1})\right\|^2\\
\texttt{MM-SVRG} & \frac{1}{bn}\sum_{i=1}^n\left\|\nabla f_i(x^k) - \nabla f_i(\tilde x^{k-1})\right\|^2\\
\texttt{MM-SARAH} & \left\|\tnablasarah f(x^{k-1}) - \nabla f(x^{k-1})\right\|^2 
\\[+.05in]
\hline
\end{array}
\end{equation}

In our analysis, we estimate $\E_k\left\|\tnabla f(x^k) - \nabla f(x^k)\right\|^2$, where $\tnabla f(x^k)$ denotes either $\tnablasaga, \tnablasvrg$, or $\tnablasarah$, and $\E_k$ denotes the conditional expectation given the history of the sequence up to the iteration $x_k$. We also provide several properties of the sequence $\{\Upsilon^k\}$. Our arguments are similar to the arguments in the study \cite{driggs2021stochastic} of stochastic proximal alternating linearized minimization algorithms using \texttt{SAGA} and \texttt{SARAH} for two-block optimization. However, we consider one block, which yields tighter bounds when compared to the bounds in \cite[Proposition 2]{driggs2021stochastic}.

\begin{lemma}\label{est:saga-svrg}
Let $\{x^k\}$ be generated by \texttt{MM-SAGA} or by \texttt{MM-SVRG} and let $\Upsilon^k$ be defined by \eqref{e:Upsilon}. Then 
\begin{enumerate}[label={\rm(\alph*)}]
    \item $\E_k\left\|\tnabla f(x^k) - \nabla f(x^k)\right\|^2 \le \Upsilon^k$.
    \item $\E_k \Upsilon^{k+1} \leq (1 - \rho) \Upsilon^k + V_\Upsilon\E_k\left\|x^{k+1} - x^k\right\|^2$, where $\rho$ and $V_\Upsilon$ are defined by \eqref{constants:saga-svrg}. 

\begin{equation}\label{constants:saga-svrg}
\arraycolsep 10pt
\def\arraystretch{1.5}
\begin{array}{l|l|l}
\hline
{\rm Method} & \rho & V_\Upsilon\\
\hline
\textup{\texttt{MM-SAGA}} & \frac{b}{2n} & \frac{(2n-b)L^2}{b^2} \\
\textup{\texttt{MM-SVRG}} & \frac{1}{2m} & \frac{(2m-1)L^2}{b}\\[+.05in]
\hline
\end{array}
\end{equation}

\end{enumerate}
\end{lemma}

\begin{proof}
First, we estimate $\E_k\left\|\tnablasaga f(x^k) - \nabla f(x^k)\right\|^2$ as follows. From the definition of the stochastic gradient estimator $\tnablasaga$ \eqref{e:tnablasaga}, it follows that
 \begin{equation}\label{est:saga-svrg:equ1}
     \begin{aligned}
         \E_k\left\|\tnablasaga f(x^k) - \nabla f(x^k)\right\|^2  & = \E_k\left\|\frac{1}{b}\sum_{i\in I_k}[\nabla f_i(x^k) - \nabla f_i(x_i^{k-1})] + v^k - \nabla f(x^k)\right\|^2\\
        & = \E_k\left\|\frac{1}{b}\sum_{i\in I_k}[\nabla f_i(x^k) - \nabla f_i(x_i^{k-1})] - [\nabla f(x^k) -v^k]\right\|^2.
     \end{aligned}
 \end{equation}
We note that for each $i\in I_k$,
\begin{equation}
\mbE_k\left[\nabla f_i(x^k) - \nabla f_i(x_i^{k-1})\right]= \nabla f(x^k)-v^k.
\end{equation}
Consequently, Lemma~\ref{l:2023feb01} implies that
\begin{align*}
\mbE_k\left\|\tnablasaga f(x^k)-\nabla f(x^k)\right\|^2
&=\frac{1}{bn}\sum_{i=1}^{n}\left\|\nabla f_i(x^k)-\nabla f_i(x_i^{k-1})\right\|^2-
\frac{1}{b}\left\|v^k-\nabla f(x^k)\right\|^2\\
&\leq \frac{1}{bn}\sum_{i=1}^{n}\left\|\nabla f_i(x^k)-\nabla f_i(x_i^{k-1})\right\|^2 = \Upsilon^k.
\end{align*}
Thus, (a) holds for \texttt{MM-SAGA}. Next, we prove that (b) holds for \texttt{MM-SAGA}. By employing the inequality
\begin{equation}\label{e:230201a-ineq}
\|x-z\|^2 \leq (1+\delta)\|x-y\|^2 + (1+\delta^{-1})\|y-z\|^2,
\end{equation}
we see that 
\begin{align*}
        \E_k\Upsilon^{k+1} =&\, \frac{1}{bn}\E_k\sum_{i=1}^n\left\|\nabla f_i(x^{k+1}) - \nabla f_i(x_i^k)\right\|^2\\
    \leq &\,\frac{1+\delta}{bn}\E_k\sum_{i=1}^n\left\|\nabla f_i(x^k) - \nabla f_i(x_i^k)\right\|^2
    + \frac{1+\delta^{-1}}{bn}\E_k\sum_{i=1}^n\left\|\nabla f_i(x^{k+1}) - \nabla f_i(x^k)\right\|^2\\
    =&\, \frac{1+\delta}{bn}\sum_{i=1}^n\left\|\nabla f_i(x^k) - \nabla f_i(x_i^{k-1})\right\|^2 - \frac{1+\delta}{bn}\E_k\sum_{i\in I_k}\left\|\nabla f_i(x^k) - \nabla f_i(x_i^{k-1})\right\|^2 \\ 
    &\,+ \frac{1+\delta^{-1}}{bn} \E_k\sum_{i=1}^n\left\|\nabla f_i(x^{k+1}) - \nabla f_i(x^k)\right\|^2
\end{align*}
\begin{align*}
&= \frac{1+\delta}{bn}\left(1-\frac{b}{n}\right)\sum_{i=1}^n\left\|\nabla f_i(x^k) - \nabla f_i(x_i^{k-1})\right\|^2  
    + \frac{1+\delta^{-1}}{bn}\E_k\sum_{i=1}^n\|\nabla f_i(x^{k+1}) - \nabla f_i(x^k)\|^2 \numberthis \label{e:2304031}\\
    &\leq (1+\delta)\left(1-\frac{b}{n}\right)\Upsilon^k  
    +  \frac{(1+\delta^{-1})L^2}{b}\E_k\left\|x^{k+1} - x^k\right\|^2 \numberthis\label{i:2304031}
\end{align*}
where 
\eqref{e:2304031} follows from
\[
\E_k\sum_{i\in I_k}\left\|\nabla f_i(x^k) - \nabla f_i(x_i^{k-1}))\right\|^2  = \frac{b}{n}\sum_{i=1}^n\left\|\nabla f_i(x^k) - \nabla f_i(x_i^{k-1})\right\|^2,
\]
and \eqref{i:2304031} follows from Assumption \ref{basic:assumption} (b). By setting $\delta = \frac{b}{2(n-b)}$, we see that
\begin{equation*}
    \E_k\Upsilon^{k+1} \leq \left(1-\frac{b}{2n}\right)\Upsilon^k + \frac{(2n-b)L^2}{b^2}\E_k\left\|x^{k+1} - x^k\right\|^2,
\end{equation*}
i.e., (b) holds for \texttt{MM-SAGA}

Next, we show that (a) holds for \texttt{MM-SVRG}. 
From the definition of $\tnablasvrg f(x^k)$ and $\tx^k$ (see \eqref{update-svrg}), it follows that
\begin{align*}
\mbE_k\left\|\tnablasvrg f(x^k)-\nabla f(x^k)\right\|^2
=\Big(1-\frac{1}{m}\Big)
\mbE_k\left\|
\frac{1}{b}\sum_{i\in I_k}\left[\nabla f_i(x^k)-\nabla f_i(\tx^{k-1})\right]
+\nabla f(\tx^{k-1})-\nabla f(x^k)\right\|^2.
\end{align*}
We note that for each $i\in I_k$,
\begin{equation}
\mbE_k\left[\nabla f_i(x^k)-\nabla f_i(\tx^{k-1})\right]
=\nabla f(x^k)-\nabla f(\tx^{k-1}).
\end{equation}
Consequently, Lemma~\ref{l:2023feb01} implies that
\begin{align*}
\mbE_k\left\|\tnablasvrg f(x^k)-\nabla f(x^k)\right\|^2
&=\Big(1-\frac{1}{m}\Big)
\frac{1}{bn}\sum_{i=1}^{n}\left\|\nabla f_i(x^k)-\nabla f_i( \tx^{k-1})\right\|^2\\
&\quad - \Big(1-\frac{1}{m}\Big)\frac{1}{b}\left\|\nabla f(x^k)-\nabla f(\tx^{k-1})\right\|^2\\
&\leq 
\Big(1-\frac{1}{m}\Big)
\frac{1}{bn}\sum_{i=1}^{n}\left\|\nabla f_i(x^k)-\nabla f_i( \tx^{k-1})\right\|^2
= \Big(1-\frac{1}{m}\Big)\Upsilon^k.
\end{align*}

Finally, we prove that (b) holds for \texttt{MM-SVRG}. Again, by recalling \eqref{e:230201a-ineq}, we see that
\begin{align*}
        \E_k\Upsilon^{k+1} &\leq \frac{1+\delta}{bn}\sum_{i=1}^n\E_k\left\|\nabla f_i(x^k) - \nabla f_i(\tilde x^k)\right\|^2 +  \frac{1+\delta^{-1}}{bn}\E_k\sum_{i=1}^n\left\|\nabla f_i(x^{k+1}) - \nabla f_i(x^k)\right\|^2\\
        & =\left(1-\frac{1}{m}\right)\frac{1+\delta}{bn}\sum_{i=1}^n\left\|\nabla f_i(x^k) - \nabla f_i(\tilde x^{k-1})\right\|^2 +  \frac{1+\delta^{-1}}{bn}\E_k\sum_{i=1}^n\left\|\nabla f_i(x^{k+1}) - \nabla f_i(x^k)\right\|^2\\
        & \leq
        \left(1-\frac{1}{m}\right)(1+\delta)\Upsilon^k + \frac{(1+\delta^{-1})L^2}{b}\E_k\left\|x^{k+1} - x^k\right\|^2.
\end{align*}
Consequently, by setting $\delta = \frac{1}{2(m-1)}$, we arrive at 
\[
\E_k\Upsilon^{k+1} \leq \left(1-\frac{1}{2m}\right)\Upsilon^k + \frac{(2m-1)L^2}{b}\E_k\left\|x^{k+1} - x^k\right\|^2,
\]
as asserted.
\end{proof}

\begin{lemma}\label{est:sarah}
Let $\{x^k\}$ be generated by \texttt{MM-SARAH} and let $\Upsilon^k$ be defined by \eqref{e:Upsilon}. Then 
\begin{enumerate}[label={\rm(\alph*)}]
    \item $\E_k\left\|\tnablasarah f(x^k) - \nabla f(x^k)\right\|^2 \le \Upsilon^k + \frac{L^2}{b}\left\|x^k - x^{k-1}\right\|^2$.
    \item $\E_k \Upsilon^{k+1} \leq (1 - \rho) \Upsilon^k + V_\Upsilon\left\|x^k - x^{k-1}\right\|^2$, where 
\end{enumerate}
\begin{equation}\label{eq:s2}
    \rho = \frac{1}{m}\quad \textrm{and}\quad V_\Upsilon = \frac{(m-1)L^2}{mb}. 
\end{equation}
\end{lemma}
\begin{proof}
(a): By recalling the definition of the \texttt{SARAH}  estimator \eqref{update-sarah}, we see that
    \begin{align*}
\E_k&\left\|\tnablasarah f(x^k) - \nabla f(x^k)\right\|^2\\
        & = \left(1 - \frac{1}{m}\right)\E_{k}\Big\|\frac{1}{b}\sum_{i\in I_k}[\nabla f_i(x^k) - \nabla f_i(x^{k-1})] + \tnablasarah f(x^{k-1}) - \nabla f(x^k)\Big\|^2\\
        & = \left(1 - \frac{1}{m}\right)\E_{k}\Big\|\frac{1}{b}\sum_{i\in I_k}\left[\nabla f_i(x^k) - \nabla f_i(x^{k-1})\right] - \nabla f(x^k) + \nabla f(x^{k-1})\\
        & \qquad\qquad\qquad\qquad\qquad + \tnablasarah f(x^{k-1}) - \nabla f(x^{k-1})\Big\|^2\\
            & =\left(1 - \frac{1}{m}\right)\E_{k}\Big\|\frac{1}{b}\sum_{i\in I_k}\left[\nabla f_i(x^k) - \nabla f_i(x^{k-1})\right] - \nabla f(x^k) + \nabla f(x^{k-1})\Big\|^2\\ 
        & \quad + \left(1 - \frac{1}{m}\right)\left\|\tnablasarah f(x^{k-1}) - \nabla f(x^{k-1})\right\|^2,
        \numberthis \label{e:2304041}
    \end{align*}
where \eqref{e:2304041} follows from the fact that for each $i\in I_k$,
\[
\E_{k}\left[\nabla f_i(x^k) - \nabla f_i(x^{k-1})\right] = \nabla f(x^k) - \nabla f(x^{k-1}).
\]
Consequently, Lemma~\ref{l:2023feb01} and the definition of $\Upsilon^k$ imply that
\begin{align}
    \E_k\left\|\tnablasarah f(x^k) - \nabla f(x^k)\right\|^2
         & = \left(1 - \frac{1}{m}\right)\frac{1}{bn}\sum_{i=1}^n\left\|\nabla f_i(x^k) - \nabla f_i(x^{k-1})\right\|^2\nonumber\\
         & \quad - \left(1 - \frac{1}{m}\right)\frac{1}{b}\left\|\nabla f(x^k) - \nabla f(x^{k-1})\right\|^2 + \left(1 - \frac{1}{m}\right)\Upsilon^k\nonumber\\
        & \leq \left(1 - \frac{1}{m}\right)\Upsilon^k + \left(1 - \frac{1}{m}\right)\frac{L^2}{b}\left\|x^k-x^{k-1}\right\|^2,\label{eq:s1}
\end{align}
where~\eqref{eq:s1} follows from Assumption \ref{basic:assumption} (b). The proof is now complete.
\end{proof}

We now present our main convergence results. The following theorem provides almost surely subsequential convergence of the iterations generated by our algorithms to a stationary point. For simplicity, we set the constant $V$ in~\eqref{constant:V}.
\begin{equation}\label{constant:V}
\arraycolsep 10pt
\def\arraystretch{1.5}
\begin{array}{l|l}
\hline
{\rm Method} &  V\\
\hline
\textup{\texttt{MM-SAGA}}/\textup{\texttt{MM-SVRG}} & 0\\
\textup{\texttt{MM-SARAH}} & \frac{L^2}{b}\\
\hline
\end{array}
\end{equation}

\begin{theorem}[Almost surely subsequential convergence]\label{t:ASC}
Let $\{x^k\}$ be a sequence generated by one of the \texttt{SVRMM} algorithms. Suppose that Assumptions~\ref{basic:assumption} and~\ref{assump:u} are satisfied, and that the condition
\begin{equation}\label{condition}
   (2\mu - L)^2 - 4(V+V_\Upsilon/\rho) > 0 
\end{equation}
where $V_
\Upsilon$ and $\rho$ are determined by~\eqref{constants:saga-svrg} for \texttt{MM-SAGA} and \texttt{MM-SVRG}, by~\eqref{eq:s2} for \texttt{MM-SARAH}, and $V$ is defined by \eqref{constant:V}, holds. Then 
\begin{enumerate}[label={\rm(\alph*)}]
    \item\label{t:ASC-a} The sequence $\{F(x^k)\}$ converges almost surely.
    \item\label{t:ASC-b} The sequence $\{\|x^k - x^{k-1}\|^2\}$ has a finite sum (in particular, vanishes) almost surely.
    \item\label{t:ASC-c} Every limit point of $\{x^k\}$ is a stationary point of $F$, almost surely.
\end{enumerate}
\end{theorem}
\begin{proof}
First, we prove \ref{t:ASC-a} and \ref{t:ASC-b}. From the definition of a surrogate function (see Definition \ref{d:surrogate}), it follows that
\begin{equation}\label{e:ASC-1}
r(x^{k+1})\leq u(x^{k+1},x^k).
\end{equation}
By combining the $L$-smoothness of $f$ with Lemma \ref{convexity-smoothness}(d), we arrive at
\begin{equation}\label{e:ASC-2}
f(x^{k+1}) \leq f(x^k) + \scal{ \nabla f(x^k)}{x^{k+1}-x^k} + \frac{L}{2}\|x^{k+1}- x^k\|^2.
\end{equation}
Now, by combining the update rule for $x^{k+1}$ with Lemma \ref{convexity-smoothness}(c), we see that
\begin{equation}\label{e:ASC-3}
\la \tnabla f(x^k), x^{k+1}-x^k\ra + u(x^{k+1},x^k) + \mu\|x^{k+1} - x^k\|^2\leq
u(x^k,x^k)=r(x^k).
\end{equation}
Consequently, by summing up \eqref{e:ASC-1}, \eqref{e:ASC-2}, and \eqref{e:ASC-3}, and by recalling that $F=f+r$, we obtain
    \begin{align*}
        F(x^{k+1}) + \frac{2\mu - L}{2}\left\|x^{k+1}-x^k\right\|^2 &\leq F(x^k) + \la \nabla f(x^k) - \tnabla f(x^k), x^{k+1} - x^k\ra\\
        & \leq F(x^k) + \frac{\eta}{2}\left\|\tnabla f(x^k) - \nabla f(x^k)\right\|^2 + \frac{1}{2\eta}\left\|x^{k+1} - x^k\right\|^2. \numberthis \label{thm1:x01}
    \end{align*}
By taking the expectation in \eqref{thm1:x01}, conditioned on $\mathcal F_k$, we arrive at
\begin{equation}\label{equ:thm:1}
    \E_k\left[F(x^{k+1}) + \left(\frac{2\mu - L}{2}-\frac{1}{2\eta}\right)\left\|x^{k+1}-x^k\right\|^2\right] \leq F(x^k) + \frac{\eta}{2}\E_k\left\|\tnabla f(x^k) - \nabla f(x^k)\right\|^2.
\end{equation}
We now consider two cases. One case, in which the sequence $\{x^k\}$ is generated  by either \texttt{MM-SAGA} or \texttt{MM-SVRG}, and the second case, in which the sequence $\{x^k\}$ is generated  by \texttt{MM-SARAH}.

{\bf Case~1}: Suppose that $\{x^k\}$ is generated  by either \texttt{MM-SAGA} or \texttt{MM-SVRG}. Lemma \ref{est:saga-svrg} (a) asserts that
\[
\E_k\left\|\tnabla f(x^k) - \nabla f(x^k)\right\|^2
   \leq \Upsilon^k + V\E_k\left\|x^{k+1}-x^k\right\|^2,
\]
which, when combined with inequality~\eqref{equ:thm:1}, implies that
\[
\begin{split}
        \E_k\left[F(x^{k+1}) + \left(\frac{2\mu - L}{2}-\frac{1}{2\eta} - \frac{\eta V}{2}\right)\left\|x^{k+1}-x^k\right\|^2\right] 
   \leq  F(x^k) + \frac{\eta}{2}\Upsilon^k.
    \end{split}
\]
By recalling Lemma \ref{est:saga-svrg} (b), which asserts that $\Upsilon^k \leq \left(\Upsilon^k - \E_k\Upsilon^{k+1}\right)/\rho + V_\Upsilon/\rho\E_k\left\|x^{k+1}-x^k\right\|^2$, we obtain 
\[
\begin{split}
        \E_k\left[F(x^{k+1}) + \left(\frac{2\mu - L}{2}-\frac{1}{2\eta} - \frac{\eta V}{2} - \frac{\eta V_\Upsilon}{2\rho}\right)\left\|x^{k+1}-x^k\right\|^2 + \frac{\eta}{2\rho}\Upsilon^{k+1}\right] 
   \leq  F(x^k)
   + \frac{\eta}{2\rho}\Upsilon^k.
    \end{split}
\]
Consequently, by setting $\Phi^{k+1} = F(x^{k+1}) + \left(\frac{2\mu - L}{2}-\frac{1}{2\eta} - \frac{\eta V}{2} - \frac{\eta V_\Upsilon}{2\rho}\right)\left\|x^{k+1}-x^k\right\|^2 + \frac{\eta}{2\rho}\Upsilon^{k+1}$, we see that
\[
\E_k \Phi^{k+1} \leq \Phi^k  - \left(\frac{2\mu - L}{2}-\frac{1}{2\eta} - \frac{\eta V}{2} - \frac{\eta V_\Upsilon}{2\rho}\right)\left\|x^k - x^{k-1}\right\|^2.
\]

{\bf Case~2:} Suppose that $\{x^k\}$ is generated  by \texttt{MM-SARAH}. It follows from \eqref{equ:thm:1} and Lemma~\ref{est:sarah} that
\begin{align*}
\E_k\left[F(x^{k+1}) + \left(\frac{2\mu - L}{2}-\frac{1}{2\eta}\right)\left\|x^{k+1}-x^k\right\|^2 + \frac{\eta}{2\rho}\Upsilon^{k+1}\right] \\
   \leq  F(x^k)
   + \frac{\eta}{2\rho}\Upsilon^k + \left(\frac{\eta V}{2} + \frac{\eta V_\Upsilon}{2\rho}\right)\left\|x^k - x^{k-1}\right\|^2.
\end{align*}
Thus, by setting $\Phi^{k+1} = F(x^{k+1}) + \left(\frac{2\mu - L}{2}-\frac{1}{2\eta}\right)\left\|x^{k+1}-x^k\right\|^2 + \frac{\eta}{2\rho}\Upsilon^{k+1}$, we obtain
\begin{equation}\label{equ:thm1:6}
    \E_k \Phi^{k+1} \leq \Phi^k  - \left(\frac{2\mu - L}{2}-\frac{1}{2\eta} - \frac{\eta V}{2} - \frac{\eta V_\Upsilon}{2\rho}\right)\left\|x^k - x^{k-1}\right\|^2.
\end{equation}
Consequently, in both cases, ~\eqref{equ:thm1:6} holds for the corresponding sequences $\{\Phi^k\}$. Now, since $$(2\mu - L)^2 - 4(V+V_\Upsilon/\rho) > 0,$$ by setting $\eta=\frac{2\mu - L}{2(V+V_\Upsilon/\rho)}$, we see that
\[
\frac{2\mu - L}{2}-\frac{1}{2\eta} - \frac{\eta V}{2} - \frac{\eta V_\Upsilon}{2\rho} = \frac{2\mu-L}{4} - \frac{V+V_\Upsilon/\rho}{2\mu-L}>0.
\] 
On the other hand, without the loss of generality, we may assume that $F^*\geq 0$. Thus, by supermartingale convergence (Lemma~\ref{supermartingale}), the sequence $\{\|x^k - x^{k-1}\|^2\}$ almost surely has a finite sum (in particular, it vanishes), and $\{\Phi^k\}$ almost surely converges to a nonnegative random variable $\Phi^{\infty}$. Consequently, by combining Lemma~\ref{est:saga-svrg}, \ref{est:sarah}, and the supermartingale convergence in Lemma~\ref{supermartingale}, $\Upsilon^k$ has a finite sum (in particular, it vanishes) almost surely. It follows that the sequence $\{F(x^k)\}$ converges to  $\Phi^{\infty}$ almost surely, which concludes the proof of \ref{t:ASC-a} and \ref{t:ASC-b}.

Proof of \ref{t:ASC-c}: 
First, we claim that 
\begin{equation}\label{e:ASC-claim}
\lim_{k\to +\infty}\left[\tnabla f(x^k) - \nabla f(x^k)\right] = 0
\quad\text{almost surely}.
\end{equation}
Indeed, by taking the total expectation in~\eqref{equ:thm1:6}, we see that
\begin{equation}\label{eq:s3}
\E\Phi^{k+1}\leq \E\Phi^k - \kappa\E\left\|x^k - x^{k-1}\right\|^2,
\end{equation}
where $\kappa = \frac{(2\mu-L)^2 - (V+V_\Upsilon/\rho)}{4(2\mu-L)}>0$. By telescoping~\eqref{eq:s3} over $k$, we arrive at
\[
\sum_{k=1}^K \kappa\E\left\|x^k - x^{k-1}\right\|^2 \leq \E\Phi^1 - \E\Phi^{K+1} \leq \E\Phi^1 - F^*,
\]
where we employed the fact that $\Phi^{K+1}\geq F(x^{K+1}) \geq F^*$. Consequently, $\left\{\E\left\|x^k - x^{k-1}\right\|^2\right\}$ has a finite sum. 

We now prove that the claim holds by considering two cases: Case~1, where $\{x^k\}$ is generated by either \texttt{MM-SAGA} or \texttt{MM-SVRG}, and case~2, where the sequence $\{x^k\}$ is generated by \texttt{MM-SARAH}. 

{\bf Case~1:} Suppose that $\{x^k\}$ is generated by either \texttt{MM-SAGA} or \texttt{MM-SVRG}. Lemma~\ref{est:saga-svrg} implies that
\begin{equation}\label{eq:s4}
\begin{split}
    \E\left\|\tnabla f(x^k) - \nabla f(x^k)\right\|^2 \le \E\Upsilon^k
    \leq \left(\E\Upsilon^k - \E\Upsilon^{k+1}\right)/\rho + V_\Upsilon/\rho\E\left\|x^{k+1}-x^k\right\|^2,
\end{split}
\end{equation}
By telescoping~\eqref{eq:s4}, we see that
\begin{align*}
    \sum_{k=0}^K\E\left\|\tnabla f(x^k) - \nabla f(x^k)\right\|^2 & \leq \left(\E\Upsilon^0 - \E\Upsilon^{K+1}\right)/\rho + V_\Upsilon/\rho\sum_{k=0}^K\E\left\|x^{k+1}-x^k\right\|^2\\
    & \leq  V_\Upsilon/\rho\sum_{k=0}^K\E\left\|x^{k+1}-x^k\right\|^2 \numberthis \label{i:2304041}
\end{align*}
where \eqref{i:2304041} is a consequence of $\Upsilon^k\geq 0$ and $\Upsilon^0 = 0$. The claim now follows from \eqref{i:2304041} and the fact that $\left\{\E\left\|x^k - x^{k-1}\right\|^2\right\}$ has a finite sum. 

{\bf Case~2:} Suppose that $\{x^k\}$ is generated by \texttt{MM-SARAH}.
Lemma~\ref{est:sarah} implies that
\begin{align*}
    \E\left\|\tnabla f(x^k) - \nabla f(x^k)\right\|^2 
    &\leq \E\Upsilon^k + V\E\left\|x^k - x^{k-1}\right\|^2\\
    &\leq \left(\E\Upsilon^k - \E\Upsilon^{k+1}\right)/\rho + \left(V+V_\Upsilon/\rho\right)\E\left\|x^k-x^{k-1}\right\|^2
    \numberthis \label{e:ASC-230504a}
\end{align*}
By telescoping \eqref{e:ASC-230504a}, we see that
\begin{align*}
    \sum_{k=0}^K\E\left\|\tnabla f(x^k) - \nabla f(x^k)\right\|^2 & \leq \left(\E\Upsilon^0 - \E\Upsilon^{K+1}\right)/\rho + \left(V+V_\Upsilon/\rho\right)\sum_{k=0}^K\E\left\|x^k-x^{k-1}\right\|^2\\
    & \leq \left(V+V_\Upsilon/\rho\right)\sum_{k=0}^K\E\left\|x^k-x^{k-1}\right\|^2 \numberthis\label{i:2304042}
\end{align*}
where \eqref{i:2304042} is a consequence of $\Upsilon^k\geq 0$ and $\Upsilon^0 = 0$. The claim now follows from \eqref{i:2304042} and the fact that $\left\{\E\left\|x^k - x^{k-1}\right\|^2\right\}$ has a finite sum.

This concludes the proof of claim \eqref{e:ASC-claim}.
 
By combining \eqref{e:ASC-claim} and part \ref{t:ASC-b}, for any sequence $x^k$ generated by one of our \texttt{SVRMM} algorithms,
\begin{equation}\label{e:ASC-4}
\lim_{k\to\infty}\tnabla f(x^k) - \nabla f(x^k) = 0\text{~~and~~}
\lim_{k\to\infty}x^k - x^{k-1} = 0
\end{equation}
almost surely.
Let $\{x^k\}$ be generated by an \texttt{SVRMM} algorithm which satisfies \eqref{e:ASC-4}.
Let $x^*$ be a limit point of $\{x^k\}$, that is, there is a subsequence $\{x^{k_j}\}$ of $\{x^k\}$ such that $x^{k_j}\to x^*$ as $j\to +\infty$. From the update rule of the \texttt{SVRMM} algorithms, it follows that
\[
\nu^{k_j+1} := -\mu(x^{k_j+1}-x^{k_j}) -  \tnabla f(x^{k_j}) \in \partial u(\cdot,x^{k_j})(x^{k_j+1}),
\]
which implies that for any $x\in\mathbb R^d$,
\begin{equation}\label{equ:thm1:2}
    u(x,x^{k_j}) \geq u(x^{k_j+1},x^{k_j}) + \la \nu^{k_j+1}, x - x^{k_j+1}\ra.
\end{equation}
By plugging $x = x^*$ into \eqref{equ:thm1:2} and by letting $j\to+\infty$, we arrive at
\begin{equation}\label{equ:thm1:3}
    r(x^*) \geq \limsup_{j\to+\infty} u(x^{k_j+1},x^{k_j}),
\end{equation}
where we employed the continuity of $u(x,y)$ in $y$ (Assumption \ref{assump:u} (a)), the fact that 
\begin{align*}
&\lim_{j\to+\infty}x^{k_j+1} = \lim_{j\to+\infty}x^{k_j} = x^*,
\text{~~and~~}\\
&\lim_{j\to+\infty}\nu^{k_j+1} = \lim_{j\to+\infty}-\tnabla f(x^{k_j}) = - \nabla f(x^*).
\end{align*}
By combining \eqref{equ:thm1:3} with the lower
semicontinuity of $u(x,y)$ in $x$ (Assumption \ref{assump:u} (b)), we conclude that
\[
\lim_{j\to+\infty} u(x^{k_j+1},x^{k_j}) = r(x^*).
\]
Consequently, by letting $j\to+\infty$ in \eqref{equ:thm1:2}, we see that for all $x\in\mathbb R^d$,
\begin{equation}\label{equ:thm1:40}
    r(x^*) \leq u(x,x^*) +  \la\nabla f(x^*), x - x^*\ra.
\end{equation}
On the other hand, since $f$ is $L$-smooth,
\begin{equation}\label{equ:thm1:4}
    f(x^*) \leq f(x) - \la\nabla f(x^*), x - x^*\ra + \frac{L}{2}\|x-x^*\|^2.
\end{equation}
By summing up \eqref{equ:thm1:40} and \eqref{equ:thm1:4}, we arrive at
\[
\begin{split}
        F(x^*) &\leq u(x,x^*) + f(x) + \frac{L}{2}\|x-x^*\|^2\\
         & = F(x) + u(x,x^*) - r(x) + \frac{L}{2}\|x-x^*\|^2\\
         & \leq F(x) + \bar h(x,x^*) + \frac{L}{2}\|x-x^*\|^2,
    \end{split}
\]
for some function $\bar h$ which satisfies Assumption \ref{assump:u} (c). Consequently, $x^*$ is a minimizer of
\begin{equation*}
    \min_{x\in \mathbb R^d} F(x) + \bar h(x,x^*) + \frac{L}{2}\|x-x^*\|^2.
\end{equation*}
It follows that
\begin{equation}\label{equ:thm1:5}
0\in \partial F(x^*) + \nabla\bar{h}(\cdot,x^*)(x^*)
=\partial F(x^*),
\end{equation}
which concludes part (c) and completes the proof.

\end{proof}
\begin{remark}[feasibility of the batchsize $b$ and the stepsize $\frac{1}{\mu}$]
For \texttt{MM-SAGA}, since $V = 0$, $V_\Upsilon = \frac{(2n-b)L^2}{b^2}$, and $\rho = \frac{b}{2n}$, condition \eqref{condition} is satisfied when
\[
(2\mu-L)^2 \geq \frac{16n^2L^2}{b^3}.
\]
For \texttt{MM-SVRG}, since $V = 0$, $V_\Upsilon = \frac{(2m-1)L^2}{b}$, and $\rho = \frac{1}{2m}$, condition \eqref{condition} is satisfied when
\[
(2\mu-L)^2 \geq \frac{16m^2L^2}{b}.
\]
For \texttt{MM-SARAH}, since $V = \frac{L^2}{b}$, $V_\Upsilon = \frac{(m-1)L^2}{mb}$, and $\rho = \frac{1}{m}$, condition \eqref{condition} is satisfied when
\[
 (2\mu-L)^2 \geq \frac{4mL^2}{b}.
\]
\end{remark}

We now provide iteration complexity in order to
obtain an $\epsilon$-stationary point. To this end, we incorporate the following additional
assumption \cite[Assumption 3(ii)]{hien2023inertial} regarding the surrogate function $u$ of $r$. 
\begin{assumption}\label{assump:unibound} For any bounded subset $\Omega$ of $\RR^d$, there exists a constant $L_u$ such that for any $x,y\in\Omega$ and for any $\nu\in \partial u(\cdot,y)(x)$, there exists $\zeta\in\partial r(x)$ such that $\|\nu-\zeta\|\leq L_u\|x - y\|$.
\end{assumption}

\begin{remark} We consider the case where $h(x,y) := u(x,y) - r(x)$ is $L_u$-smooth in $x$ and $\nabla h(\cdot,y)(y)=0$, as assumed in~\cite{Mairal2015}. We assert that Assumption \ref{assump:unibound} captures this case. Indeed, since $\partial u(\cdot,y)(x) = \partial r(x) +\nabla h(\cdot,y)(x)$, if $\nu\in \partial u(\cdot,y)(x)$, then there exists $\zeta\in \partial r(x)$ such that
\[
 \|\nu - \zeta\| =  \left\|\nabla h(\cdot,y)(x)\right\| = \|\nabla h(\cdot,y)(x) - \nabla h(\cdot,y)(y)\| \leq L_u\|x - y\|.
\]
It follows that the proximal surrogates (Examples~\ref{ex:prox-surr}), the Lipschitz gradient surrogates (Examples~\ref{ex:lip-grad-surr}), and the DC surrogates (Example~\ref{ex:dc-surr}) satisfy Assumption \ref{assump:unibound}. Furthermore, Assumption \ref{assump:unibound} captures composite surrogates (Example~\ref{ex:comp_surr}) with non-smooth approximation error functions. Indeed, let $\nu\in \partial u(\cdot,y)(x)$. Then 
$$\nu = \big( \eta'_1(g_1(y_1))\xi_1,\ldots, \eta'_m(g_m(y_m))\xi_m\big)\quad \text{where}\quad \xi_j\in\partial g_j(x_j).$$ 
On the other hand, it follows from \cite[Corollary 5Q]{Rockafellar1981TheTO} that 
\[
\partial (\eta_j\circ g_j)(x_j) = \eta'_j(g_j(x_j))\partial g_j(x_j).
\]
Consequently, by letting $\zeta = \big(\eta'_1(g_j(x_1))\xi_1,...,\eta'_m(g_m(x_m))\xi_m\big)\in\partial r(x)$, we arrive at
\[
\|\nu - \zeta\| \leq \sqrt{\sum_{j=1}^mL_\eta^2\|g_j(x_j) - g_j(y_j)\|^2\|\xi_j\|^2} \leq L_\eta L_g^2\|x -y\|
\]
which implies that Assumption \ref{assump:unibound} is satisfied by letting $L_u = L_\eta L_g^2$.
\end{remark}

\begin{theorem}[Iteration complexity]
Let $\{x^k\}$ be a sequence generated by one of the \texttt{SVRMM} algorithms. Suppose that Assumptions \ref{basic:assumption}, \ref{assump:u}, \ref{assump:unibound}, and condition \eqref{condition} are satisfied. Then for any $K\in\mathbb{N}$,
\begin{equation}
    \min_{k = 1,2,..,K}\E\dist^2(0,\partial F(x^k)) \leq \frac{8(2\mu-L)[(L+\mu+L_u)^2+V_\Upsilon/\rho](F(x^0) - F^*)}{K[(2\mu-L)^2  - 4(V+V_\Upsilon/\rho)]} = \mathcal O(1/K).
\end{equation}

In other words, it takes at most $\mathcal O(1/\epsilon^2)$ iterations, in expectation, to obtain an $\epsilon$-stationary point of $F$. 
\end{theorem}
\begin{proof}
From the update rule of the \texttt{SVRMM} algorithms, it follows that
\[
\nu^{k+1} := -\mu(x^{k+1}-x^{k}) -  \tnabla f(x^{k}) \in \partial u(\cdot,x^{k})(x^{k+1}).
\]
Thus, by Assumption \ref{assump:unibound}, there exists $\zeta^{k+1}\in\partial r(x^{k+1})$ such that $\|\nu^{k+1} - \zeta^{k+1}\| \leq L_u\|x^{k+1}-x^{k}\|$. Consequently, by invoking  Lemma \ref{convexity-smoothness} (b), it follows that
\[
\nabla f(x^{k+1}) + \zeta^{k+1} \in \nabla f(x^{k+1}) + \partial r(x^{k+1}) = \partial F(x^{k+1}).
\]
We see that
\[
\begin{split}
        \dist(0,\partial F(x^{k+1})) &\leq \left\|\nabla f(x^{k+1}) + \zeta^{k+1}\right\|\\ 
        & = \left\|\nabla f(x^{k+1}) - \nabla f(x^k) + \nabla f(x^k) - \tnabla f(x^k) - \mu(x^{k+1}-x^k) + \zeta^{k+1} - \nu^{k+1}\right\| \\ 
        &\leq \left(L+\mu+L_u\right)\left\|x^{k+1}-x^k\right\| + \left\|\tnabla f(x^k) - \nabla f(x^k)\right\|.
    \end{split}
\]
It now follows that 
\begin{equation}\label{equ:thm2:1}
    \begin{split}
        \E\dist^2\left(0,\partial F(x^{k+1})\right)  
        \leq 2\left(L+\mu+L_u\right)^2\E\left\|x^{k+1}-x^k\right\|^2 + 2\E\left\|\tnabla f(x^k) - \nabla f(x^k)\right\|^2.
    \end{split}
\end{equation}
On the other hand, if the sequence is generated by either \texttt{MM-SAGA} or \texttt{MM-SVRG}, Lemma~\ref{est:saga-svrg} implies that
\begin{equation}\label{equ:thm2:1-1}
\E\left\|\tnabla f(x^k) - \nabla f(x^k)\right\|^2 \leq \frac{1}{\rho}\left(\E\Upsilon^k - \E\Upsilon^{k+1}\right) + \frac{V_\Upsilon}{\rho}\E\left\|x^{k+1}-x^k\right\|^2.
\end{equation}
By plugging~\eqref{equ:thm2:1-1} into \eqref{equ:thm2:1} we arrive at
\begin{equation}\label{equ:thm2:3}
    \E\dist^2\left(0,\partial F(x^{k+1})\right)  
        \leq \left[2\left(L+\mu+L_u\right)^2+\frac{2V_\Upsilon}{\rho}\right]\E\left\|x^{k+1}-x^k\right\|^2 +\frac{2}{\rho}\left(\E\Upsilon^k - \E\Upsilon^{k+1}\right).
\end{equation}
Moreover, it follows from \eqref{equ:thm1:6} that
\begin{equation}\label{equ:boundxk}
    \E\left\|x^{k+1}-x^k\right\|^2\leq \kappa\left(\E\Phi^{k+1} - \E\Phi^{k+2}\right),
\end{equation}
where $\kappa = \frac{4(2\mu-L)}{(2\mu-L)^2 - 4(V+V_\Upsilon/\rho)}$. By combining \eqref{equ:thm2:3} and \eqref{equ:boundxk}, we see that
\begin{equation}\label{eq:s5}
\E\dist^2\left(0,\partial F(x^{k+1})\right)  
        \leq 2\kappa\left[\left(L+\mu+L_u\right)^2+V_\Upsilon/\rho\right](\E\Phi^{k+1} - \E\Phi^{k+2})  + \frac{2}{\rho}\left(\E\Upsilon^k - \E\Upsilon^{k+1}\right).
\end{equation}
Consequently, by telescoping~\eqref{eq:s5} over $k$, we obtain
\begin{align*}
\sum_{k=1}^K\E\dist^2\left(0,\partial F(x^k)\right)  
        & \leq  2\kappa\left[\left(L+\mu+L_u\right)^2+V_\Upsilon/\rho\right]\left(\E\Phi^1 - \E\Phi^{K+1}\right) + \frac{2}{\rho}\left(\E\Upsilon^0 - \E\Upsilon^K\right)\\
        & \leq 2\kappa\left[\left(L+\mu+L_u\right)^2+V_\Upsilon/\rho\right](\Phi^1 - F^*), \numberthis \label{i:2304043}
    \end{align*}
where \eqref{i:2304043} follows from $\Upsilon^0 = 0$, $\Upsilon^k \geq 0$, and $\Phi^k\geq F(x^k)\geq F^*$. We conclude the proof in the case of \texttt{MM-SAGA} and \texttt{MM-SVRG} by combining~\eqref{i:2304043} and the fact that
\begin{align*}
    \Phi^1 & = F(x^1) + \left(\frac{2\mu-L}{4} -  \frac{V+V_\Upsilon/\rho}{2\mu-L}\right)\left\|x^1 - x^0\right\|^2 + \frac{2\mu-L}{4(\rho V+V_\Upsilon)}\Upsilon^1\\
    & \leq F(x^1) + \frac{2\mu-L}{4}\left\|x^1 - x^0\right\|^2 + \frac{2\mu-L}{4(\rho V+V_\Upsilon)}\Upsilon^1\\
    & \leq F(x^1) + \frac{2\mu-L}{2}\left\|x^1 - x^0\right\|^2 \numberthis\label{i:2304045}\\
    & \leq F(x^0), \numberthis\label{i:2304044}
\end{align*}
where \eqref{i:2304044} follows from \eqref{thm1:x01} and \eqref{i:2304045} follows by recalling that
\[
\frac{1}{\rho V+V_\Upsilon} \Upsilon^1= \frac{1}{\rho V+V_\Upsilon} \frac{1}{bn}\sum_{i=1}^n\left\|\nabla f_i(x^1) - \nabla f_i(x^0)\right\|^2
    \leq \frac{1}{\rho V+V_\Upsilon} \frac{L^2}{b}\left\|x^1-x^0\right\|^2\leq \left\|x^1-x^0\right\|^2.
\]

In the case where the sequence $\{x^k\}$ is generated by \texttt{MM-SARAH}, Lemma~\ref{est:sarah} implies that 
\begin{equation}\label{equ:thm2:1-2}
\E\left\|\tnablasarah f(x^k) - \nabla f(x^k)\right\|^2 = \E\Upsilon^{k+1} \leq \frac{1}{\rho}\left(\E\Upsilon^{k+1} - \E\Upsilon^{k+2}\right) + \frac{V_\Upsilon}{\rho}\E\left\|x^{k+1}-x^k\right\|^2.
\end{equation}
By plugging~\eqref{equ:thm2:1-2} into \eqref{equ:thm2:1}, we see that
\begin{equation}\label{eq:s6}
\E\dist^2(0,\partial F(x^{k+1}))  
        \leq \left[2\left(L+\mu+L_u\right)^2+\frac{2V_\Upsilon}{\rho}\right]\E\left\|x^{k+1}-x^k\right\|^2 +\frac{2}{\rho}\left(\E\Upsilon^{k+1} - \E\Upsilon^{k+2}\right).
\end{equation}
By telescoping~\eqref{eq:s6} over $k$, we obtain
\[
\begin{split}
        \sum_{k=1}^K\E\dist^2\left(0,\partial F(x^k)\right)  
        & \leq 2\left[\left(L+\mu+L_u\right)^2+V_\Upsilon/\rho\right]\sum_{k=0}^{K-1}\E\left\|x^{k+1}-x^k\right\|^2 + \frac{2}{\rho}\left(\E\Upsilon^1 - \E\Upsilon^{K+1}\right)\\
        & \leq  2\kappa\left[\left(L+\mu+L_u\right)^2+V_\Upsilon/\rho\right]\left(\E\Phi^1 - \E\Phi^{K+1}\right)\\
        & \leq 2\kappa\left[\left(L+\mu+L_u\right)^2+V_\Upsilon/\rho\right]\left(\Phi^1 - F^*\right),
    \end{split}
\]
where we employed the fact that $\Upsilon^1 = 0$, $\Upsilon^k \geq 0$, $\Phi^k\geq F(x^k)\geq F^*$, and \eqref{equ:boundxk}, which concludes the proof in the case of \texttt{MM-SARAH}. 
\end{proof}

The \texttt{SVRMM} algorithm includes two parameters: $\mu$ and the batch-size $b$. Setting $\mu = L$ yields a larger stepsize $\frac{1}{L}$ compared to other stochastic gradient methods such as \texttt{DCA-SAGA} and \texttt{DCA-SVRG} \cite{le2022stochastic}, with a stepsize of $\frac{1}{2L}$, \texttt{ProxSVRG} \cite{j2016proximal}, with a stepsize of $\frac{1}{3L}$, and \texttt{ProxSVRG+} \cite{li2018simple}, with a stepsize of $\frac{1}{6L}$. After fixing $\mu$, we select the batch-size that satisfies condition \eqref{condition}. The following corollary summarizes choices of the batch size $b$ and the associated complexity of each algorithm in terms of the number of individual stochastic gradient valuations $\nabla f_i$.
\begin{corollary}\label{cor:complexity}
\begin{enumerate}[label={\rm(\alph*)}]
    \item In the case of \texttt{MM-SAGA}, we set $\mu = L$ and $b = \lfloor 4^{2/3}n^{2/3}\rfloor$. Consequently, the complexity  is $\mathcal O(n^{2/3}/\epsilon^2)$.
    \item In the case of \texttt{MM-SVRG}, we set $\mu = L$, $b = \lfloor n^{2/3}\rfloor$ and $m =  \frac{\sqrt{b}}{4} =  \frac{n^{1/3}}{4}$. Consequently, the complexity is $\mathcal O(n^{2/3}/\epsilon^2)$.
    \item In the case of \texttt{MM-SARAH}, we set $\mu = L$, $b = \lfloor n^{1/2}\rfloor$ and $m =  \frac{b}{4} = \frac{n^{1/2}}{4}$. Consequently, the complexity is $\mathcal O(n^{1/2}/\epsilon^2)$.
\end{enumerate}
\end{corollary}

\begin{remark}
\begin{enumerate}[label={\rm(\alph*)}]
    \item Our results coincide with the best-known complexity bounds to obtain an $\epsilon$-stationary point in expectation for \texttt{ProxSAGA} \cite{j2016proximal}, \texttt{ProxSVRG} \cite{j2016proximal,li2018simple}, \texttt{SPIDER} \cite{fang2018spider}, \texttt{SpiderBoost} \cite{wang2019spiderboost}, and \texttt{ProxSARAH} \cite{pham2020proxsarah} methods in the particular case where $r = 0$ \cite{fang2018spider,wang2019spiderboost} or the case where $r$ is convex \cite{j2016proximal,li2018simple,pham2020proxsarah}.
    \item 
    Within Corollary \ref{cor:complexity}, instead of setting $\mu=L$, we may first fix a batchsize $b$, then choose a compatible parameter $\mu$ to comply with condition \eqref{condition}. In particular, we may pick a batch-size $b\in\{1,2,...,n-1\}$, then set $\mu = (4nL/b^{3/2} + L)/2$ for \texttt{MM-SAGA}, $\mu = (4mL/b^{1/2} + L)/2$ for \texttt{MM-SVRG} and $\mu = (2m^{1/2}L/b^{1/2} + L)/2$ for \texttt{MM-SARAH}. 
\end{enumerate}

\end{remark}

\section{Numerical experiments}\label{numerical_exp}
We now examine the applicability and efficiency of our \texttt{SVRMM} algorithms. To this end, we consider the following three problems: sparse binary classification with nonconvex loss and regularizer, sparse multi-class logistic regression with nonconvex regularizer, and feedforward neural network training. 

We compare six algorithms:
\begin{itemize}
    \item \texttt{MM-SAGA} with $\mu = L$ and $b = \lfloor 4^{2/3}n^{2/3}\rfloor$;
    \item \texttt{MM-SVRG} with $\mu = L$, $m = \frac{n^{1/3}}{4}$, and $b = \lfloor n^{2/3}\rfloor$;
    \item \texttt{MM-SARAH} with $\mu = L$, $m = \frac{n^{1/2}}{4}$ and $b = \lfloor n^{1/2}\rfloor$;
    \item \texttt{SDCA} \cite{le2020stochastic} with $\mu = 1.1L$ and  $b=\lfloor n/10\rfloor$, which performed well in \cite{le2020stochastic};
    \item \texttt{DCA-SAGA} \cite{le2022stochastic} with $\mu = 2L$ and $b = \lfloor 2\sqrt{n\sqrt{n+1}}\rfloor$; 
    \item \texttt{DCA-SVRG} \cite{le2022stochastic} with $\mu = 2L$, $b = \lfloor n^{2/3}\rfloor$, and the inner loop length $M = \lfloor \frac{\sqrt{b}}{4\sqrt{e-1}}\rfloor$.
\end{itemize}

In our experiments, we run each algorithm 20 epochs repeated 20 times, where each epoch consists of $n$ gradient evaluations. We are interested
in the relative loss residuals $\frac{F(w^k) - F^*}{|F^*|}$, where $F^*$ is the minimum loss values generated by all algorithms, and in classification accuracy on testing sets. 

All tests are performed using Python on a Linux server with configuration: Intel(R) Xeon(R) Gold 5220R CPU 2.20GHz
of 64GB RAM. The code is available at \url{https://github.com/nhatpd/SVRMM}.

\subsection{Sparse binary classification with nonconvex loss and regularizer}
Let $\{(a_i,b_i): i = 1,...,n\}$ be a training set with observation vectors $a_i\in \mathbb{R}^d$ and labels $b_i \in \{-1,1\}$. We consider the sparse binary classification with nonconvex loss function and nonconvex regularizer:
\begin{equation}\label{binaryclassification}
    \min_{w\in\mathbb R^d}\biggl\{F(w) = \frac{1}{n}\sum_{i=1}^n\ell(a_i^Tw,b_i) + r(w) \biggr\},
\end{equation}
where $\ell$ is a nonconvex loss function and $r$ is a regularization term.  We revisit a nonconvex loss function from \cite{zhao2010convex}: $\ell(s,t) = \biggl(1 - \frac{1}{1+\exp(-ts)}\biggr)^2$, and the exponential regularization from~\cite{Bradley1998}: $r(w) = \sum_{i=1}^d\eta\circ g(w_i)$, where $\eta$ and $g$ are the functions
\begin{equation}\label{exp:penaty}
    \eta(t) = \lambda\biggl(1 -\exp(-\alpha t)\biggr) \quad \text{and}\quad g(w) = |w|,
\end{equation}
where $\lambda$ and $\alpha$ are nonnegative tuning parameters. The hessian matrix of $\ell(a_i^Tw,b_i)$ is evaluated as follows:
\[
\nabla^2\ell(a_i^Tw,b_i) = \frac{4\exp(2b_ia_i^Tw) - 2\exp(b_ia_i^Tw)}{\left(\exp(2b_ia_i^Tw) + 1\right)^4}a_ia_i^T.
\]
We thus have
\[
\|\nabla^2\ell(a_i^Tw,b_i)\| = \frac{\left|4\exp(2b_ia_i^Tw) - 2\exp(b_ia_i^Tw)\right|}{\left(\exp(2b_ia_i^Tw) + 1\right)^4}\|a_ia_i^T\| \leq  \frac{39+55\sqrt{33}}{2304}\|a_i\|^2.
\]
Therefore, $\ell(a_i^Tw,b_i)$ is $L$-smooth with
\[
L = \frac{39+55\sqrt{33}}{2304}\max\limits_{i=1,...,n}\|a_i\|^2,
\]
and, in this case, problem \eqref{binaryclassification} is within the scope of problem \eqref{model} when we let $f_i(w) = \ell(a_i^Tw,b_i)$. Moreover, since $\eta$ is concave and $\lambda\alpha^2$-smooth on $\mathbb R_+$, and since $g$ is convex and $1$-Lipschitz continuous, we set a surrogate function $u$ for $r$ as follows: 
\[
u(w,w^k) = r(w^k)+ \sum_{i=1}^d\lambda\alpha\exp(-\alpha|w_i^k|)(|w_i|-|w_i^k|).
\]
Assumptions \ref{assump:u} and \ref{assump:unibound} are then satisfied. The \texttt{SVRMM} algorithms update $w^{k+1}$ to be the solution of the nonsmooth convex subproblem:
\[
\min_{w\in\mathbb R^d}\frac{\mu}{2}\|w-w^k\|^2 + \la \tnabla f(w^k),w\ra +  \sum_{i=1}^d\lambda\alpha\exp(-\alpha|w_i^k|)|w_i|,
\]
for which a closed-form solution was provided in \cite[Section 6.5.2]{parikh2014proximal} by
\[
 w^{k+1}_i = \max\biggl\{|v_i^k|,\lambda\alpha\exp(-\alpha|w_i^k|)/\mu\biggr\}\text{sign}(v_i^k),
\]
where $v_i^k = w_i^k - \tnabla f(x^k)/\mu$. 

\subsection{Sparse multi-class logistic regression with nonconvex regularizer}
We revisit the multi-class logistic regression with a nonconvex regularizer:
\begin{equation}\label{multi-class}
    \min_{W\in\mathbb R^{d\times q}}\biggl\{F(W) = \frac{1}{n}\sum_{i=1}^n\ell(b_i,a_i,W) + r(W)\biggr\},
\end{equation}
where $q$ is the number of classes, $\{(a_i,b_i): i = 1,2,...,q\}$ is a training set with the feature vectors $a_i\in\mathbb R^d$ and the labels $b_i\in\{1,2,...,q\}$,  $r$ is a regularizer, and $\ell(b_i,a_i,W)$ is a loss function defined by
\[
\ell(b_i,a_i,W) =\log\biggl(\sum_{k=1}^q\exp(a_i^Tw_k)\biggr) -a_i^Tw_{b_i},
\]
where $w_k$ is the $k$-th column of $W$. We employ an exponential $\ell_2$ regularizer, defined by
\[
r(W) = \lambda\sum_{i=1}^d\eta(\|W_i\|),
\]
where $\eta$ is defined as in \eqref{exp:penaty}, and $W_i$ is the $i$-th row  of $W$. Since $\ell(b_i,a_i,W)$ is $L$-smooth with $L=\frac{q-1}{q}\max_{i=1,..,n}\|a_i\|^2$, it follows that problem \eqref{multi-class} is within the scope of problem \eqref{model} when we set $f_i(W) = \ell(b_i,a_i,W)$. The \texttt{SVRMM} algorithms applied to \eqref{multi-class} iteratively determine a surrogate function of $r(W)$ at $W^k$ by
\[
 u(W,W^k) = r(W^k) + \sum_{i=1}^d\lambda\alpha\exp(-\alpha\|W_i^k\|)(\|W_i\| - \|W_i^k\|),
\]
and then update $W^{k+1}$ by
\[
W^{k+1} = \arg\min_{W\in\mathbb R^{d\times q}} \frac{\mu}{2}\|W - W^k\|^2 + \la\tnabla f(W^k),W\ra + \sum_{i=1}^d\lambda\alpha\exp(-\alpha\|W_i^k\|)\|W_i\|
\]
for which a closed-form solution was provided in \cite[Section 6.5.1]{parikh2014proximal} by
\[
W_i^{k+1} = \begin{cases} \biggl(1 - \frac{\lambda\alpha\exp(-\alpha\|W_i^k\|)/\mu}{\|V_i^k\|}\biggr)V_i^k \quad& \text{if}\quad  \|V_i^k\| \geq \lambda\alpha\exp(-\alpha\|W_i^k\|)/\mu,\\
    0 & \text{otherwise}, \end{cases}
\]
where $V_i^k = W^k_i - \tnabla f(W^k)_i/\mu$.

\subsection{Feedforward neural network training problem with nonconvex regularizer}
We consider the nonconvex optimization model arising in a feedforward neural network configuration
\begin{equation}\label{nn}
    \min_{w\in\mathbb R^D}\biggl\{F(w) = \frac{1}{n}\sum_{i=1}^n\ell(h(w,a_i),b_i) + r(w) \biggr\},
\end{equation}
where all of the weight matrices and bias vectors of the neural network are concatenated in one vector of variables $w$, $(a_i,b_i)_{i=1}^n$ is a training data set with the feature vectors $a_i\in\mathbb R^d$ and the labels $b_i\in\{1,2,...,q\}$,  $h$ is a composition of linear transforms and activation functions of the form $h(w,a) = \sigma_l(W_l\sigma_{l-1}(W_{l-1}\sigma_{l-2}(\cdots \sigma_0(W_0a + c_0)\cdots)+c_{l-1}) + c_l)$, where $W_i$ is a weight matrix, $c_i$ is a bias vector, $\sigma_i$ is an activation function, $l$ is the number of layers, $\ell$ is the soft-max cross-entropy loss, and $r$ is a regularizer. By considering the exponential regularization $r(w) = \sum_{i=1}^D\eta\circ g(w_i)$, where $\eta$ and $g$ are set in \eqref{exp:penaty}, problem \eqref{nn} is within the scope of of problem \eqref{model} when we let $f_i(w) = \ell(h(w,a_i),b_i)$. The \texttt{SVRMM} algorithms applied to \eqref{nn} are different from the \texttt{SVRMM} algorithms for problem \eqref{binaryclassification} only in computation of stochastic gradient estimates $\tnabla f$. In our experiment, we employ a one-hidden-layer fully connected neural network, $784\times 100\times 10$, as studied in \cite{pham2020proxsarah}. The activation function $\sigma_i$ of the hidden layer is ReLU. 

\subsection{Experiment setups and data sets}
In our experiments, for the first two problems \eqref{binaryclassification} and \eqref{multi-class}, all of the algorithms under study start at the zero point, while for the last problem \eqref{nn}, we use the global\_variables\_initializer function from Tensorflow. We set the regularization parameters $\alpha = 5$ for the first two problems and $\alpha = 0.05$ for the latter, and we fix $\lambda = 1/n$. These regularization parameters are standard in the literature, e.g., \cite{pham2020proxsarah,Bradley1998}. It is important to mention that in all of the experiments, we use the same problem settings for all of the algorithms. 

For the sparse binary classification, we carried out the experiments
on six different well-known data sets: a9a, w8a, rcv1, real-sim, epsilon, and url. For the sparse multi-class logistic regression, we test all of the algorithms on four data sets: dna, shuttle, Sensorless, and connect-4. Finally, for the feedforward neural network training, we use two data sets minist and fashion\_mnist, to compare \texttt{MM-SVRG} and \texttt{MM-SARAH} with \texttt{DCA-SVRG}. We randomly pick 90\% of the data for training and the rest for testing. The characteristics of the data sets
are provided in Table \ref{datasets}. The first eleven data sets are obtained from the LIBSVM Data website\footnote{\url{https://www.csie.ntu.edu.tw/~cjlin/libsvmtools/datasets/}} while the data sets minist and fashion\_mnist are obtained from the library tensorflow.keras.datasets. 
\begin{table}[]
\centering
\caption{Data sets used in experiments}\label{datasets}
\begin{tabular}{llll}
\hline\noalign{\smallskip}   Data set & Data points & Class number &Feature number\\
\hline\noalign{\smallskip}
a9a 		& 32,561 	& 2 & 123\\
w8a & 49,749 & 2 & 300\\
rcv1-binary  	& 20,242 	& 2 & 47,236\\
real-sim & 27,309 & 2 & 20,958\\
epsilon & 400,000 & 2 & 2,000\\
url & 2,396,130 & 2 & 3,231,961\\
\hline
dna & 2,000 & 11 & 180\\ 
shuttle & 43,500 & 7 & 9\\
Sensorless & 58,509 & 11 & 48\\
connect-4 & 67,557 & 3 & 126\\
mnist & 60,000 & 10 & 780\\
fashion mnist & 60,000 & 10 & 780\\
\hline
\end{tabular}
\end{table}
\subsection{Results}
\begin{figure*}[!htpb]
\begin{center}
\begin{tabular}{cc}
\includegraphics[width=0.42\linewidth]{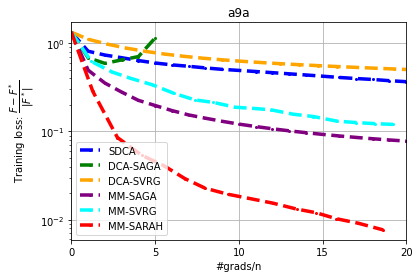} & 
\includegraphics[width=0.42\linewidth]{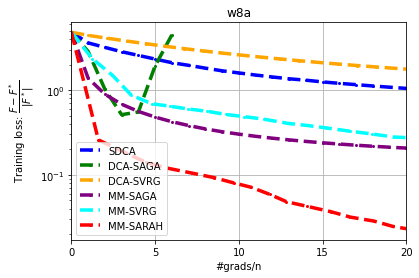} \\ \includegraphics[width=0.42\linewidth]{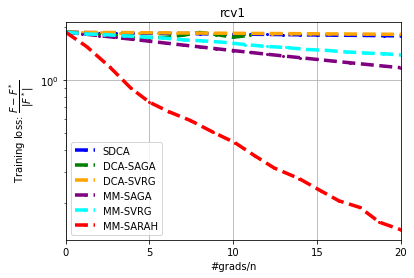} &
\includegraphics[width=0.42\textwidth]{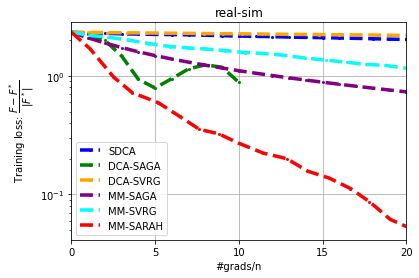}  \\
\includegraphics[width=0.42\linewidth]{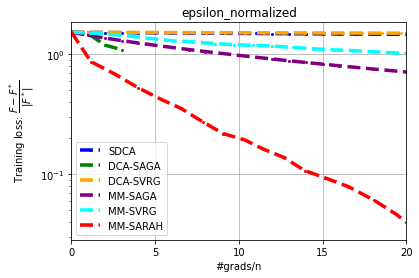} & \includegraphics[width=0.42\linewidth]{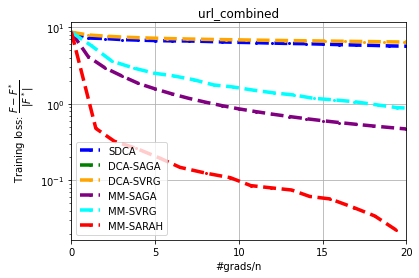}\\
\includegraphics[width=0.42\textwidth]{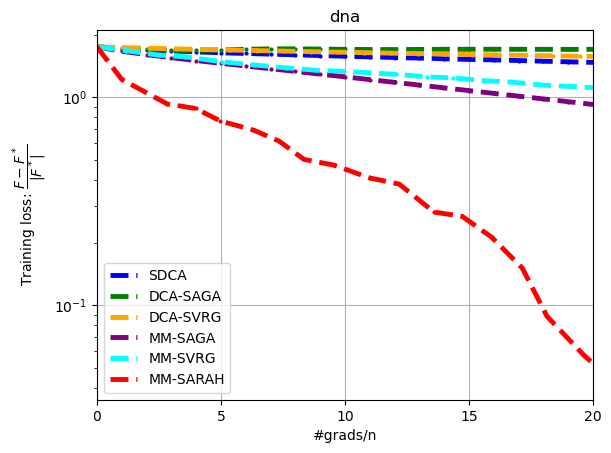} & 
\includegraphics[width=0.42\textwidth]{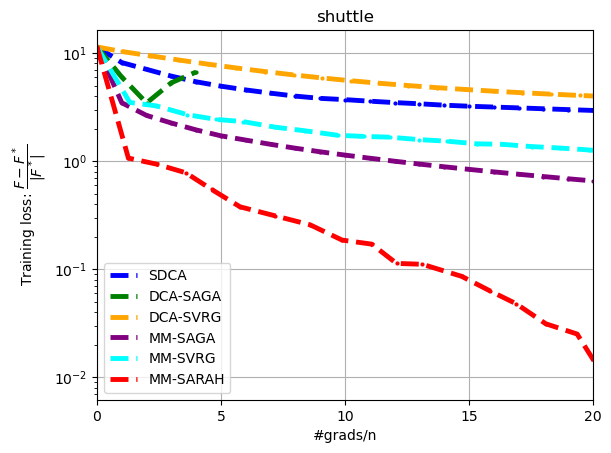} \\ \includegraphics[width=0.43\linewidth]{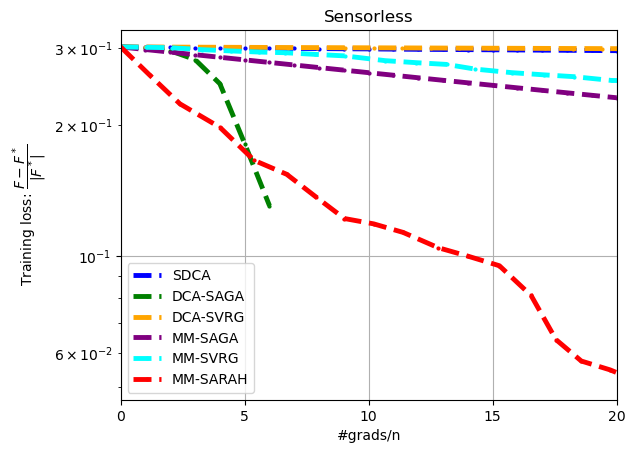} &
\includegraphics[width=0.42\linewidth]{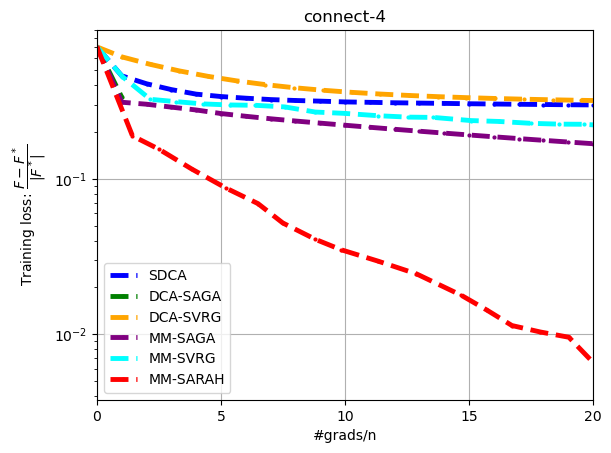} 
\end{tabular}
\caption{Evolution of the average value of the relative loss residuals with respect to the epoch on a9a, w8a, rcv1, real-sim, epsilon\_normalized, url\_combined, dna, shuttle, Sensorless, and connect-4.
\label{fig:results}} 
\end{center}
\vspace{-0.2in}
\end{figure*}

\begin{figure*}[!htpb]
\begin{center}
\begin{tabular}{cc}
\includegraphics[width=0.42\linewidth]{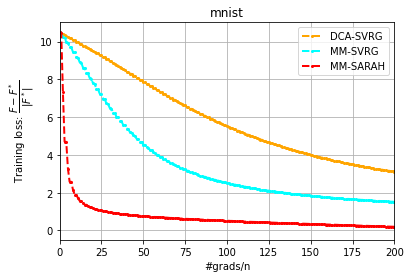} 
\includegraphics[width=0.42\linewidth]{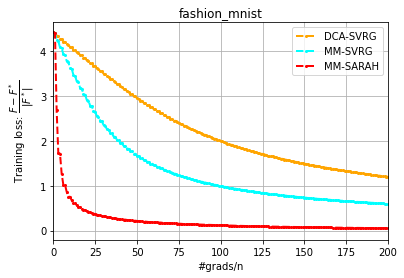} 

\end{tabular}
\caption{Evolution of the average value of the relative loss residuals with respect to the epoch on mnist and fashion mnist.
\label{fig:mnist}} 
\end{center}
\vspace{-0.2in}
\end{figure*}

We plotted the curves of the average value of relative loss residuals versus epochs in Figures \ref{fig:results}-\ref{fig:mnist}, and reported the average and the standard deviation
of the relative loss residuals and the testing accuracy in Table \ref{table:performance}. We observe from Table \ref{table:performance} and  Figures \ref{fig:results}-\ref{fig:mnist} that \texttt{MM-SARAH} has the fastest convergence on all of the data sets. \texttt{MM-SARAH} achieves not only the best relative loss residuals but also the best classification accuracy on the testing sets. This illustrates the theoretical results, see Corollary \ref{cor:complexity}, where \texttt{MM-SARAH} has the best complexity among these algorithms. In addition, \texttt{MM-SAGA} performs better than \texttt{DCA-SAGA}, which is not stable on the first ten data sets. This illustrates the benefit of the proximal term in the iterate of \texttt{MM-SAGA}. Moreover, \texttt{MM-SVRG} performs better than \texttt{DCA-SVRG} on all of the data sets, which illustrates the benefit of the loop-less variant of \texttt{SVRG} in \texttt{MM-SVRG}.

\section{Conclusion}\label{conclusion}
We introduced three stochastic variance-reduced MM algorithms: \texttt{MM-SAGA}, \texttt{MM-SVRG}, and \texttt{MM-SARAH}, combining the MM principle and the variance reduction techniques from \texttt{SAGA}, \texttt{SVRG}, and \texttt{SARAH} for solving a class of nonconvex nonsmooth optimization problems with the large-sum structure. The complex objective function is approximated by compatible surrogate functions, providing closed-form solutions in the updates of our algorithms. At the same time, we employ the benefits of the stochastic gradient estimators (\texttt{SAGA}, loop-less \texttt{SVRG}, and loop-less \texttt{SARAH}) to overcome the challenge of the large-sum structure. We provided almost surely subsequential convergence of \texttt{MM-SAGA}, \texttt{MM-SVRG}, and \texttt{MM-SARAH} to a stationary point under mild assumptions. In addition, we proved that our algorithms possess the state-of-the-art complexity bounds in terms of the number of gradient evaluations without assuming that the approximation errors of the regularizer $r$ are $L$-smooth. We applied our new algorithms to three important problems in machine learning in order to demonstrate the advantages of combining the \texttt{MM} principle with \texttt{SAGA}, \texttt{SVRG}, and \texttt{SARAH}. Overall, \texttt{MM-SARAH} outperforms other stochastic algorithms under consideration. This is not surprising since the methods based on \texttt{SAGA} and \texttt{SVRG} have unvoidable limitations. In particular, \texttt{SAGA} requires storing the most recent gradient of each component function $f_i$ while \texttt{SVRG} employs a pivot iterate $\tilde{x}^k$ that may be unchanged during many iterations and, thus, may no longer  be highly correlated with the current iterate~$x^k$. 

\begin{table}[!htpb]
    \centering
\resizebox{0.99\textwidth}{!}{
    \begin{tabular}{l|llllll} 
\toprule 
\multirow{ 2}{*}{Data} & \multicolumn{6}{c}{Method}\\
\cline{2-7}
 & \texttt{SDCA} & \texttt{DCA-SAGA} & \texttt{DCA-SVRG} & \texttt{MM-SAGA} & \texttt{MM-SVRG} & \texttt{MM-SARAH} \\
 \midrule
 \multirow{2}{*}{w7a} & 0.962 (0.063) & 0.584 (0.390) & 1.580 (0.019) & 0.208 (0.054) & 0.276 (0.034) & \textbf{0.047} (0.029)\\
 &  0.893 (0.007) & 0.893 (0.007) & \textbf{0.895} (0.005) & 0.893 (0.007) & \textbf{0.895} (0.005) & 0.894 (0.005)\\
 \midrule
 \multirow{2}{*}{a9a} &  0.366 (0.012) & 0.591 (0.063) & 0.5 (0.002) & 0.078 (0.016) & 0.12  (0.012) & \textbf{0.008} (0.004)\\
 &  0.779 (0.007) & 0.782 (0.032) & 0.758 (0.005) & 0.833 (0.009) & 0.834 (0.005) & \textbf{0.845} (0.006)\\
 \midrule
 \multirow{2}{*}{w8a} &  1.048 (0.053) & 0.513 (0.38) & 1.742 (0.008) & 0.207 (0.05) & 0.273 (0.069) & \textbf{0.022} (0.015)\\
 &  0.895 (0.006) & 0.895 (0.006) & 0.89 (0.004) & \textbf{0.896} (0.006) & 0.89  (0.004) & 0.892 (0.003)\\
 \midrule
 \multirow{2}{*}{rcv1} & 1.78 (0.025) & 1.74 (0.091) & 1.816 (0.005) & 1.173 (0.153) & 1.342 (0.111) & \textbf{0.14} (0.108)\\
 &  0.854 (0.043) & 0.856 (0.044) & 0.877 (0.031) & 0.891 (0.029) & 0.901 (0.02) & \textbf{0.927} (0.011) \\
 \midrule
 \multirow{2}{*}{real-sim} & 2.034 (0.162) & 0.791 (0.142) & 2.198 (0.106) & 0.734 (0.158) & 1.118 (0.085) & \textbf{0.05} (0.026)\\
 &  0.797 (0.025) & 0.872 (0.014) & 0.784 (0.005) & 0.882 (0.012) & 0.858 (0.011) & \textbf{0.927} (0.004)\\
 \midrule
 \multirow{2}{*}{epsilon} & 1.468 (0.076) & 1.095 (0.046) & 1.495 (0.071) & 0.712 (0.051) & 1.003 (0.061) & \textbf{0.034} (0.025)\\
 &  0.71 (0.007) & 0.713 (0.004) & 0.707 (0.003) & 0.844 (0.009) & 0.8 (0.013) & \textbf{0.885} (0.002)\\
 \midrule
 \multirow{2}{*}{url} & 5.644 (0.262) & 8.584 (1.751) & 6.266 (0.201) & 0.469 (0.184) & 0.866 (0.15) & \textbf{0.023} (0.014)\\
 &  0.721 (0.015) & 0.67 (0.011) & 0.671, (0.001) & 0.962 (0.007) & 0.961 (0.001) & \textbf{0.969} (0.001)\\
 \midrule
 \multirow{2}{*}{dna}  & 1.471 (0.029) & 1.68 (0.002) & 1.56 (0.015) & 0.923 (0.032) & 1.056 (0.074) & \textbf{0.042} (0.012)\\
 &  0.514 (0.027) & 0.514 (0.027) & 0.491 (0.031) & 0.719 (0.049) & 0.631 (0.004) & \textbf{0.928} (0.015)\\
 \midrule
 \multirow{2}{*}{shuttle} & 2.942 (0.063) & 3.451 (0.095) & 3.775 (0.008) & 0.655 (0.034) & 1.162 (0.007) & \textbf{0.009} (0.005)\\
 &  0.802 (0.007) & 0.802 (0.007) & 0.797 (0.009) & 0.902 (0.006) & 0.856 (0.007) & \textbf{0.949} (0.005)\\
 \midrule
 \multirow{2}{*}{Sensorless}  & 0.296 (0.001) & 0.13 (0.134) & 0.299 (0.006) & 0.231 (0.003) & 0.249 (0.007) & \textbf{0.051} (0.029)\\
 &  0.291 (0.033) & 0.302 (0.028) & 0.267 (0.025) & 0.314 (0.032) & 0.338 (0.006) & \textbf{0.48} (0.008)\\
 \midrule
 \multirow{2}{*}{connect-4} & 0.299 (0.001) & 0.334 (0.013) & 0.314 (0.001) & 0.168 (0.005) & 0.214 (0.011) & \textbf{0.005} (0.002)\\
 &  0.659 (0.004) & 0.659 (0.004) & 0.663 (0.001) & 0.679 (0.007) & 0.664 (0.002) & \textbf{0.742} (0.003)\\
 \midrule
 \multirow{2}{*}{mnist} & - & - & 1.672 (0.245) & - & 1.021 (0.119) & \textbf{0.028} (0.013)\\
  &  - & - & 0.88 (0.009) & - & 0.905 (0.007) & \textbf{0.954} (0.002)\\
 \midrule
 \multirow{2}{*}{fashion mnist} & - & - & 0.665 (0.112) & - & 0.355 (0.065) & \textbf{0.015} (0.008)\\
 &  - & -& 0.753 (0.016) & - & 0.799 (0.01) & \textbf{0.846} (0.002)\\
 \bottomrule
\end{tabular}
}
    \caption{Performance of the comparative algorithms reported in the form of the mean value (standard deviation) with respect to the relative loss residuals (upper row) and testing accuracy (lower row) for each of the data sets. Bold values indicate the best results.}
    \label{table:performance}
    \vspace{-0.3in}
\end{table}

\newpage 
{\bf Acknowledgements.} DNP, SB, and HMP are partially supported by a Seed Grant from the Kennedy College of Sciences, University of Massachusetts Lowell. SB is partially supported by a Simons Foundation Collaboration Grant for Mathematicians. DNP and HMP are partially supported by a Gift from Autodesk, Inc. NG is partially supported by the National Science Foundation grant: NSF DMS \#2015460.


\begin{thebibliography}{99}

\bibitem{allen2018natasha}
Zeyuan Allen-Zhu.
\newblock Natasha 2: Faster non-convex optimization than sgd.
\newblock {\em Advances in neural information processing systems}, 31, 2018.

\bibitem{allen2018neon2}
Zeyuan Allen-Zhu and Yuanzhi Li.
\newblock Neon2: Finding local minima via first-order oracles.
\newblock {\em Advances in Neural Information Processing Systems}, 31, 2018.

\bibitem{attcon}
Hedy Attouch, J{\'e}r{\^o}me Bolte, and Benar~Fux Svaiter.
\newblock Convergence of descent methods for semi-algebraic and tame problems:
  proximal algorithms, forward--backward splitting, and regularized
  gauss--seidel methods.
\newblock {\em Mathematical Programming}, 137(1):91--129, Feb 2013.

\bibitem{bauschke_combettes_2017}
Heinz~H Bauschke, Patrick~L Combettes, et~al.
\newblock {\em Convex analysis and monotone operator theory in Hilbert spaces},
  volume 408.
\newblock Springer, 2017.

\bibitem{becafa}
A.~Beck and M.~Teboulle.
\newblock A fast iterative shrinkage-thresholding algorithm for linear inverse
  problems.
\newblock {\em SIAM J. Imag. Sci.}, 2:183--202, 2009.

\bibitem{beck2017first}
Amir Beck.
\newblock {\em First-order methods in optimization}.
\newblock SIAM, 2017.

\bibitem{Bradley1998}
P.~S. Bradley and O.~L. Mangasarian.
\newblock Feature selection via concave minimization and support vector
  machines.
\newblock In {\em Proceeding of international conference on machine learning
  ICML'98}, 1998.

\bibitem{Candes_reweighting}
E.~J. Cand\`es, M.~B. Wakin, and S.~P. Boyd.
\newblock Enhancing sparsity by reweighted $\ell_1$ minimization.
\newblock {\em J. Fourier Anal. Appl.}, 14(5--6):877--905, Dec. 2008.

\bibitem{chouzenoux2022sabrina}
Emilie Chouzenoux and Jean-Baptiste Fest.
\newblock Sabrina: A stochastic subspace majorization-minimization algorithm.
\newblock {\em Journal of Optimization Theory and Applications}, pages 1--34,
  2022.

\bibitem{chouzenoux2017stochastic}
Emilie Chouzenoux and Jean-Christophe Pesquet.
\newblock A stochastic majorize-minimize subspace algorithm for online
  penalized least squares estimation.
\newblock {\em IEEE Transactions on Signal Processing}, 65(18):4770--4783,
  2017.

\bibitem{combettes2011proximal}
Patrick~L Combettes and Jean-Christophe Pesquet.
\newblock Proximal splitting methods in signal processing.
\newblock In {\em Fixed-point algorithms for inverse problems in science and
  engineering}, pages 185--212. Springer, 2011.

\bibitem{defazio2014saga}
Aaron Defazio, Francis Bach, and Simon Lacoste-Julien.
\newblock Saga: A fast incremental gradient method with support for
  non-strongly convex composite objectives.
\newblock {\em Advances in neural information processing systems}, 27, 2014.

\bibitem{dempster1977maximum}
Arthur~P Dempster, Nan~M Laird, and Donald~B Rubin.
\newblock Maximum likelihood from incomplete data via the em algorithm.
\newblock {\em Journal of the Royal Statistical Society: Series B
  (Methodological)}, 39(1):1--22, 1977.

\bibitem{driggs2021stochastic}
Derek Driggs, Junqi Tang, Jingwei Liang, Mike Davies, and Carola-Bibiane
  Schonlieb.
\newblock A stochastic proximal alternating minimization for nonsmooth and
  nonconvex optimization.
\newblock {\em SIAM Journal on Imaging Sciences}, 14(4):1932--1970, 2021.

\bibitem{fang2018spider}
Cong Fang, Chris~Junchi Li, Zhouchen Lin, and Tong Zhang.
\newblock Spider: Near-optimal non-convex optimization via stochastic
  path-integrated differential estimator.
\newblock {\em Advances in Neural Information Processing Systems}, 31, 2018.

\bibitem{geman1995nonlinear}
Donald Geman and Chengda Yang.
\newblock Nonlinear image recovery with half-quadratic regularization.
\newblock {\em IEEE transactions on Image Processing}, 4(7):932--946, 1995.

\bibitem{hale2008fixed}
Elaine~T Hale, Wotao Yin, and Yin Zhang.
\newblock Fixed-point continuation for $\backslash$ell\_1-minimization:
  Methodology and convergence.
\newblock {\em SIAM Journal on Optimization}, 19(3):1107--1130, 2008.

\bibitem{hien2022inertial}
Le~Thi~Khanh Hien, Duy~Nhat Phan, and Nicolas Gillis.
\newblock Inertial alternating direction method of multipliers for non-convex
  non-smooth optimization.
\newblock {\em Computational Optimization and Applications}, 83(1):247--285,
  2022.

\bibitem{hien2023inertial}
LTK Hien, DN~Phan, and N~Gillis.
\newblock An inertial block majorization minimization framework for nonsmooth
  nonconvex optimization.
\newblock {\em Journal of Machine Learning Research}, 24:1--41, 2023.

\bibitem{j2016proximal}
Sashank J~Reddi, Suvrit Sra, Barnabas Poczos, and Alexander~J Smola.
\newblock Proximal stochastic methods for nonsmooth nonconvex finite-sum
  optimization.
\newblock {\em Advances in neural information processing systems}, 29, 2016.

\bibitem{johnson2013accelerating}
Rie Johnson and Tong Zhang.
\newblock Accelerating stochastic gradient descent using predictive variance
  reduction.
\newblock {\em Advances in neural information processing systems}, 26, 2013.

\bibitem{kappro}
A.~Kaplan and R.~Tichatschke.
\newblock Proximal point methods and nonconvex optimization.
\newblock {\em Journal of Global Optimization}, 13(4):389--406, Dec 1998.

\bibitem{khanh2022block}
Le~Thi Khanh~Hien, Duy~Nhat Phan, Nicolas Gillis, Masoud Ahookhosh, and
  Panagiotis Patrinos.
\newblock Block bregman majorization minimization with extrapolation.
\newblock {\em SIAM Journal on Mathematics of Data Science}, 4(1):1--25, 2022.

\bibitem{kovalev2020don}
Dmitry Kovalev, Samuel Horv{\'a}th, and Peter Richt{\'a}rik.
\newblock Don’t jump through hoops and remove those loops: Svrg and katyusha
  are better without the outer loop.
\newblock In {\em Algorithmic Learning Theory}, pages 451--467. PMLR, 2020.

\bibitem{lange2000optimization}
Kenneth Lange, David~R Hunter, and Ilsoon Yang.
\newblock Optimization transfer using surrogate objective functions.
\newblock {\em Journal of computational and graphical statistics}, 9(1):1--20,
  2000.

\bibitem{letthe}
H.~A Le~Thi and T.~Pham~Dinh.
\newblock The {DC} (difference of convex functions) programming and {DCA}
  revisited with {DC} models of real world nonconvex optimization problems.
\newblock {\em Annals of Operations Research}, 133:23--46, 2005.

\bibitem{le2020stochastic}
Hoai~An Le~Thi, Hoai~Minh Le, Duy~Nhat Phan, and Bach Tran.
\newblock Stochastic dca for minimizing a large sum of dc functions with
  application to multi-class logistic regression.
\newblock {\em Neural Networks}, 132:220--231, 2020.

\bibitem{le2022stochastic}
Hoai~An Le~Thi, Hoang Phuc~Hau Luu, Hoai~Minh Le, and Tao~Pham Dinh.
\newblock Stochastic dca with variance reduction and applications in machine
  learning.
\newblock {\em Journal of Machine Learning Research}, 23(206):1--44, 2022.

\bibitem{li2018simple}
Zhize Li and Jian Li.
\newblock A simple proximal stochastic gradient method for nonsmooth nonconvex
  optimization.
\newblock {\em Advances in neural information processing systems}, 31, 2018.

\bibitem{Mairal2015}
Julien Mairal.
\newblock Incremental majorization-minimization optimization with application
  to large-scale machine learning.
\newblock {\em SIAM Journal on Optimization}, 25(2):829--855, 2015.

\bibitem{marbre}
B.~Martinet.
\newblock Brève communication. r\'egularisation d'in\'equations
  variationnelles par approximations successives.
\newblock {\em ESAIM: Mathematical Modelling and Numerical
  Analysis-Mod\'elisation Math\'ematique et Analyse Num\'erique},
  4(R3):154--158, 1970.

\bibitem{neal1998view}
Radford~M Neal and Geoffrey~E Hinton.
\newblock A view of the em algorithm that justifies incremental, sparse, and
  other variants.
\newblock In {\em Learning in graphical models}, pages 355--368. Springer,
  1998.

\bibitem{nesterov2013gradient}
Yu~Nesterov.
\newblock Gradient methods for minimizing composite functions.
\newblock {\em Mathematical programming}, 140(1):125--161, 2013.

\bibitem{nesterov2003introductory}
Yurii Nesterov.
\newblock {\em Introductory lectures on convex optimization: A basic course},
  volume~87.
\newblock Springer Science \& Business Media, 2003.

\bibitem{nesterov2018lectures}
Yurii Nesterov et~al.
\newblock {\em Lectures on convex optimization}, volume 137.
\newblock Springer, 2018.

\bibitem{nguyen2017sarah}
Lam~M Nguyen, Jie Liu, Katya Scheinberg, and Martin Tak{\'a}{\v{c}}.
\newblock Sarah: A novel method for machine learning problems using stochastic
  recursive gradient.
\newblock In {\em International Conference on Machine Learning}, pages
  2613--2621. PMLR, 2017.

\bibitem{nguyen2019optimal}
Lam~M Nguyen, Marten van Dijk, Dzung~T Phan, Phuong~Ha Nguyen, Tsui-Wei Weng,
  and Jayant~R Kalagnanam.
\newblock Optimal finite-sum smooth non-convex optimization with sarah.
\newblock {\em arXiv preprint arXiv:1901.07648}, 2019.

\bibitem{parikh2014proximal}
Neal Parikh, Stephen Boyd, et~al.
\newblock Proximal algorithms.
\newblock {\em Foundations and trends{\textregistered} in Optimization},
  1(3):127--239, 2014.

\bibitem{parizi2019generalized}
Sobhan~Naderi Parizi, Kun He, Reza Aghajani, Stan Sclaroff, and Pedro
  Felzenszwalb.
\newblock Generalized majorization-minimization.
\newblock In {\em International Conference on Machine Learning}, pages
  5022--5031. PMLR, 2019.

\bibitem{pham2020proxsarah}
Nhan~H Pham, Lam~M Nguyen, Dzung~T Phan, and Quoc Tran-Dinh.
\newblock Proxsarah: An efficient algorithmic framework for stochastic
  composite nonconvex optimization.
\newblock {\em J. Mach. Learn. Res.}, 21(110):1--48, 2020.

\bibitem{phacon}
T.~Pham~Dinh and H.~A. Le~Thi.
\newblock Convex analysis approach to {D.C. programming}: Theory, algorithms
  and applications.
\newblock {\em Acta Mathematica Vietnamica}, 22(1):289--355, 1997.

\bibitem{razaviyayn2013unified}
Meisam Razaviyayn, Mingyi Hong, and Zhi-Quan Luo.
\newblock A unified convergence analysis of block successive minimization
  methods for nonsmooth optimization.
\newblock {\em SIAM Journal on Optimization}, 23(2):1126--1153, 2013.

\bibitem{robbins1951stochastic}
Herbert Robbins and Sutton Monro.
\newblock A stochastic approximation method.
\newblock {\em The annals of mathematical statistics}, pages 400--407, 1951.

\bibitem{robbins1971convergence}
Herbert Robbins and David Siegmund.
\newblock A convergence theorem for non negative almost supermartingales and
  some applications.
\newblock In {\em Optimizing methods in statistics}, pages 233--257. Elsevier,
  1971.

\bibitem{VariationalAnalysis}
R.~Rockafellar and R.~Wets.
\newblock {\em Variational Analysis}.
\newblock Springer Berlin Heidelberg, 2009.

\bibitem{rocmon}
R.~Tyrrell Rockafellar.
\newblock Monotone operators and the proximal point algorithm.
\newblock {\em SIAM Journal on Control and Optimization}, 14(5):877--898, 1976.

\bibitem{Rockafellar1981TheTO}
R.~Tyrrell Rockafellar.
\newblock {\em The theory of subgradients and its applications to problems of
  optimization : convex and nonconvex functions}.
\newblock 1981.

\bibitem{schmidt2017minimizing}
Mark Schmidt, Nicolas Le~Roux, and Francis Bach.
\newblock Minimizing finite sums with the stochastic average gradient.
\newblock {\em Mathematical Programming}, 162(1):83--112, 2017.

\bibitem{Lethi2017}
Hoai An~Le Thi, Hoai~Minh Le, Phan~Duy Nhat, and Bach Tran.
\newblock Stochastic {DCA} for the large-sum of non-convex functions problem
  and its application to group variable selection in classification.
\newblock In {\em Proceedings of the 34th International Conference on Machine
  Learning, {ICML} 2017, Sydney, NSW, Australia, 6-11 August 2017}, volume~70
  of {\em Proceedings of Machine Learning Research}, pages 3394--3403, 2017.

\bibitem{wang2019spiderboost}
Zhe Wang, Kaiyi Ji, Yi~Zhou, Yingbin Liang, and Vahid Tarokh.
\newblock Spiderboost and momentum: Faster variance reduction algorithms.
\newblock {\em Advances in Neural Information Processing Systems}, 32, 2019.

\bibitem{zhang2010nearly}
Cun-Hui Zhang.
\newblock Nearly unbiased variable selection under minimax concave penalty.
\newblock {\em The Annals of statistics}, 38(2):894--942, 2010.

\bibitem{zhao2010convex}
Lei Zhao, Musa Mammadov, and John Yearwood.
\newblock From convex to nonconvex: a loss function analysis for binary
  classification.
\newblock In {\em 2010 IEEE International Conference on Data Mining Workshops},
  pages 1281--1288. IEEE, 2010.

\end{thebibliography}
\end{document}